\documentclass[11pt]{amsart}
\usepackage{amsmath,amssymb,amsthm}
\usepackage[margin=1.1in]{geometry}
\usepackage{hyperref}
\newtheorem{theorem}{Theorem}[section]
\newtheorem{lemma}[theorem]{Lemma}
\newtheorem{proposition}[theorem]{Proposition}
\newtheorem{corollary}[theorem]{Corollary}
\theoremstyle{definition}
\newtheorem{definition}[theorem]{Definition}
\newtheorem{notation}[theorem]{Notation}
\newtheorem{example}[theorem]{Example}
\newtheorem{remark}[theorem]{Remark}
\newtheorem{problem}[theorem]{Problem}
\newcommand{\R}{\mathcal R}
\newcommand{\A}{\mathcal A}
\newcommand{\C}{\mathbb C}
\newcommand{\comm}{\operatorname{comm}}
\newcommand{\qnil}{\operatorname{qnil}}
\newcommand{\gd}{\mathrm{gd}}
\newcommand{\Dz}{\mathrm D}
\newcommand{\ind}{\operatorname{ind}}
\newcommand{\diag}{\operatorname{diag}}
\newcommand{\gap}{\operatorname{gap}}

\usepackage{xcolor}
\definecolor{linkblue}{RGB}{0,0,128} 
\hypersetup{
    colorlinks=true,
    linkcolor=linkblue,
    citecolor=linkblue,
    urlcolor=linkblue
}
\title[Spectral coincidence and generalized Hua inversion]{Spectral coincidence and generalized Hua inversion in rings and Banach algebras}
\author{Yaashaa Golovanov}
\address{AMMOC, An International Math Circle in Honor of Vladimir Igorevich Arnold and Jerrold Eldon Marsden}
\email{golovanov@ammoc.org}
\thanks{Preprint, \today.}
\date{}
\begin{document}

\begin{abstract}
For additively coupled elements with invertible sum and selected commuting
outer inverses, we study coincidence of the core idempotents and inversion
of the Hua element. Intertwining forces coincidence for finite-Drazin data
in unital rings and for generalized-Drazin data in complex unital Banach
algebras. Under coincidence, the output splits into an invertible core
and a complementary parallel sum, called the residue, yielding exact
criteria for inverse existence and success of the prescribed candidate.
Explicit ring examples show that strict generalized-Drazin input inverses
need not force coincidence under intertwining, and that a common core
with a Harte-quasinilpotent residue need not give a strict output inverse.
For rings finite as modules over a commutative central subring with the
same identity, we prove intertwining transfer and automatic bicommutation
for commuting outer inverses with Harte-quasinilpotent remainders,
without Noetherianity. This gives input-core coincidence under
intertwining. If the central subring is Noetherian, every
Harte-quasinilpotent element also has separated left and right powers.
For common-core data in the finite-central class, the candidate is a
strict output inverse exactly when the residual parallel sum is
Harte-quasinilpotent. In the Banach setting, uniform input spectral gaps characterize
continuity of the core and candidate along convergent input pairs, with
intertwining and quasinilpotent input remainders at every term and at the
limit. The output gap is equivalent when the residues of the sequence and
limit are also quasinilpotent. Exact defect identities give local estimates
on augmented inverse data. At an exactly intertwined base with invertible
coefficient operators, we determine the feasible first-order defect space
and its sharp norm constant, with a norm-specific matrix calibration.
\end{abstract}
\maketitle
\section{Introduction}\label{sec:intro}
Hua's identity has the algebraic form
\[
 a-aba=a(a+c)^{-1}c=(a^{-1}+c^{-1})^{-1},\qquad a+c=b^{-1},
\]
when the displayed inverses exist. Replacing the input inverses by generalized inverses presents two different difficulties. The inputs may carry incompatible core idempotents. Even if these idempotents coincide, the complementary part of the output need not satisfy the condition required for the proposed inverse. In the present paper, we treat these two issues separately. Here spectral coincidence means coincidence of the core idempotents, not equality of the input spectra.

We now fix the notation. Let $\R$ be a unital associative ring. For $a\in\R$ and $b\in\R^\times$, put $c=b^{-1}-a$, and choose commuting outer inverses $a',c'$:
\[
 a'aa'=a',\quad aa'=a'a=p,\qquad c'cc'=c',\quad cc'=c'c=r.
\]
The input remainders, candidate inverse, and Hua element are
\[
 n_a=a(1-p),\qquad n_c=c(1-r),\qquad s=a'+c',\qquad t=a-aba.
\]
Our intertwining condition is the element equation $pb=br$. It is not an additional commutativity assumption on $a$ and $b$.
\begin{theorem}[Coincidence and its inverse consequences]\label{thm:headline}
In the setting above, suppose either that both input remainders are nilpotent, or that $\R$ is a complex unital Banach algebra and both are quasinilpotent. Then
\[
 pb=br\quad\Longleftrightarrow\quad p=r.
\]
Under these equivalent conditions, write $e=p=r$, $f=1-e$, and $B=f\R f$. All input and output data are block diagonal along $e$, and
\[
 t=T_0+\nu,\quad T_0\in(e\R e)^\times,\quad (T_0)^{-1}_{e\R e}=s,\quad \nu=n_a(n_a+n_c)^{-1}_B n_c.
\]
We call the complementary component $\nu=ftf=t-T_0$
the \emph{residue}, or \emph{residual parallel sum}, of $t$
relative to the common core idempotent $e$.
The inverse in the displayed expression for $\nu$ is taken
in the corner ring $B=f\R f$, whose identity is $f$.
This terminology imposes no nilpotence or quasinilpotence
condition on $\nu$.

In every ring $t$ has a finite Drazin inverse if and only if $\nu$ has one in $B$, in which case $t^{\Dz}=s+\nu_B^{\Dz}$; the corresponding assertion holds for the group inverse. In the Banach case,
\[
 t\in\R^{\gd}\quad\Longleftrightarrow\quad\nu\in B^{\gd}.
\]
When these inverses exist, $t^{\gd}=s+\nu_B^{\gd}$.
In particular, $s=t^\#$ exactly when $\nu=0$, $s=t^{\Dz}$ exactly when $\nu$ is nilpotent, and, in the Banach case, $s=t^{\gd}$ exactly when $\nu$ is quasinilpotent.
\end{theorem}
The proof follows by combining Theorems~\ref{thm:rigidity},
\ref{thm:common}, \ref{thm:output}, and~\ref{thm:spectrum}.
This inverse framework connects three principal parts of the paper.

\emph{Central finiteness, separation and bicommutation.}
If $\R$ is finite as a module over a commutative central subring
with the same identity, Theorem~\ref{thm:central-transfer} transfers
intertwiners through commuting outer inverses with Harte-quasinilpotent
remainders. Its proof eliminates the off-diagonal recurrences using
integrality of the inverse coefficients; no Noetherian hypothesis is
needed. It gives input-core coincidence under intertwining and,
independently of Hua coupling, automatic bicommutation of each given
witness (Corollary~\ref{cor:central-bicommutation}). If the central
subring is Noetherian, Theorem~\ref{thm:central-separation} gives the
stronger conclusion that all left and right power intersections of a
Harte-quasinilpotent element vanish. The example following
Theorem~\ref{thm:central-transfer} distinguishes these conclusions.

\emph{Independent ring obstructions and residual behavior.}
Example~\ref{ex:ring} has strict generalized-Drazin input inverses and
intertwining, but distinct input cores. Example~\ref{rem:ring-output}
has a common nonzero proper core, finite-Drazin inputs of index three
and a Harte-quasinilpotent residual parallel sum, yet no strict output
inverse. Central module-finiteness excludes both mechanisms without
Noetherianity. It does not force residual smallness: in that class, for
common-core data the candidate is a strict output inverse precisely
when the residual parallel sum is Harte-quasinilpotent
(Proposition~\ref{prop:separated-output}).
Section~\ref{sec:examples} separates the further output regimes.

\emph{Continuity recovery and feasible perturbations.}
The recovery identities in Proposition~\ref{prop:canonical} give the
full-class continuity criterion of Theorem~\ref{thm:full-continuity},
even when the output inverse does not exist. Along convergent input
pairs with every term and the limit in the full class, either uniform
input gap characterizes continuity of $e$ and $s$; the output criterion
requires the restricted class. For approximate intertwining,
Section~\ref{sec:defects} gives exact defect identities on augmented
inverse data. At an exactly intertwined base with invertible coefficient
operators, Theorem~\ref{thm:tangent-sharp} realizes every off-diagonal
first-order direction and identifies the sharp norm constant.
The defect and tangent estimates concern augmented inverse data;
they do not assert unrestricted raw-input stability or backward-error
repair.

The Hua element $a(a+c)^{-1}c$ is a parallel sum, following the
matrix operation of Anderson and Duffin~\cite{AndersonDuffin69}.
Equivalently, it is the Schur complement of $a+c$ in
$\left(\begin{smallmatrix}a&a\\a&a+c\end{smallmatrix}\right)$.
These algebraic formulas do not carry positivity or electrical-network
hypotheses into the non-normal setting. The prescribed inverse is a Bott--Duffin inverse, with the two-idempotent ring formulation in Drazin~\cite[Definition 3.2 and Proposition 3.3]{Drazin12}. Drazin's inverse-intertwining theorem~\cite[Theorem 2.2]{Drazin13} is a classical mechanism behind the finite-index argument. The particular extra relation here is
\[
 c(p-r)=-(p-r)a,
\]
which is supplied by additive coupling and the stated intertwining. The absorption law underlying the common-core equivalence has antecedents for arbitrary outer inverses in Bapat, Jain, Karantha and Raj~\cite[Theorem~8(i),(iii)]{BapatEtAl17}, and for canonical generalized-Drazin data in Koliha and Patr\'icio~\cite[Theorem~6.1(v)]{KolihaPatricio02}. The precise specializations are recorded before Theorem~\ref{thm:common}.

Li, Chen and Mosi\'c~\cite[Theorems 2.3 and 2.5]{LiChenMosic20} prove reflexive-inverse Hua identities under range and null-space hypotheses. Their witnesses need not commute. The present outer data need not satisfy the inner equation $aa'a=a$, and their hypotheses are not identified with our intertwining condition. Under the notation of those two theorems, the relevant comparison is $T=-c$, $T^+=-c'$, $S=-s$, and $ABA-A=-t$. In the commuting-reflexive subcase the remainders vanish; with our intertwining condition, the configuration is classical. We do not classify all pairs satisfying a bare equality $t^{\Dz}=a^{\Dz}+c^{\Dz}$ outside the prescribed-core framework.

The Banach inverse calculus uses Koliha's generalized-Drazin
theory~\cite{Koliha96}. Its spectral machinery and the general continuity
tools are classical. The restricted output criterion is an application
of Koliha--Rako\v{c}evi\'c~\cite[Theorem 2.4]{KolihaRakocevic98}
through canonicity, whereas the full-class criterion uses the input
recovery identities. The coupled residual formula determines the
output-inverse question to which the spectral theory is applied.

\section{Outer data and the exact common-core structure}\label{sec:common}
Throughout the algebraic sections $\R$ is a unital associative ring. For $x\in\R$, write $\comm(x)=\{y:xy=yx\}$ and
\[
 \comm^2(x)=\{z:zy=yz\text{ for every }y\in\comm(x)\}.
\]
The core idempotent of inverse data is denoted by $p$ or $e$; its complement, not the core, is the zero-spectral idempotent in the analytic convention. The corner $e\R e$ has identity $e$. The zero corner is the zero ring, whose sole element is its identity and its inverse. In the analytic setting its spectrum is empty.
\begin{definition}[Commuting outer data]\label{def:outer}
A commuting outer inverse of $a$ is $a'$ satisfying $a'aa'=a'$ and $aa'=a'a$. Put $p=aa'$, $A_0=pa$, and $n_a=a-A_0$. We call $A_0$ the core and $n_a$ the input remainder. These words by themselves impose no nilpotence condition.
\end{definition}
\begin{lemma}[The input corner]\label{lem:input}
For these data, $p^2=p$, $a'=pa'p$, $A_0\in(p\R p)^\times$ with inverse $a'$, and
\[
 a=A_0+n_a,\qquad n_a\in(1-p)\R(1-p),\qquad a'n_a=n_aa'=A_0n_a=n_aA_0=0.
\]
\end{lemma}
\begin{proof}
The outer equation gives $p^2=a(a'aa')=aa'=p$ and $pa'=a'p=a'$. Since $a$ commutes with $p$, its two off-diagonal corners vanish. The products $A_0a'=a'A_0=p$ and the remaining annihilation identities follow by multiplying the complementary corners.
\end{proof}
A finite Drazin inverse is a commuting outer inverse whose input remainder is nilpotent. It is unique and belongs to the double commutant of its input~\cite[Theorem 1]{Drazin58}. The group inverse is the case of zero input remainder. We use index zero for units; a zero element in a nonzero ring has index one. In a general ring, $q$ is Harte-quasinilpotent when $1+qy$ is invertible for every $y\in\comm(q)$; see~\cite[Definition 2.2]{Koliha96} and~\cite[Definition 2.1]{KolihaPatricio02}. The strict ring generalized Drazin inverse additionally requires its witness to lie in $\comm^2(a)$. Mere commutation with a Harte-quasinilpotent remainder does not substitute for that requirement. In a complex Banach algebra, the commuting-outer definition with spectral-radius-zero remainder does give the unique generalized Drazin inverse~\cite[Theorems 4.2--4.3]{Koliha96}.

The ring and analytic quasinilpotence conventions agree in each nonzero complex
unital Banach algebra $\A$~\cite[Definition 2.2 and the following note]{Koliha96}:
\[
 q\in\qnil(\A)\quad\Longleftrightarrow\quad\rho_{\A}(q)=0.
\]
Indeed, if $\rho_{\A}(q)=0$ and $y\in\comm(q)$, then
$\|(qy)^k\|^{1/k}\le\|q^k\|^{1/k}\|y\|\to0$. The Neumann series, therefore,
inverts $1+qy$.

Conversely, Harte quasinilpotence applied to $y=-\lambda^{-1}1$ makes
$\lambda1-q$ invertible for every $\lambda\ne0$. The spectrum in a nonzero
complex unital Banach algebra is nonempty, so $\sigma_{\A}(q)=\{0\}$ and
$\rho_{\A}(q)=0$. The same argument applies in each nonzero Banach corner
with its own identity. In the zero corner the spectrum is empty and the sole
element is Harte-quasinilpotent; its spectral radius is taken to be zero,
consistently with the root-norm formula.

\begin{notation}[Coupled data]\label{not:coupled}
Unless otherwise stated, $b\in\R^\times$, $a+c=b^{-1}$, $a',c'$ are commuting outer inverses, and
\[
 p=aa',\quad r=cc',\quad n_a=a(1-p),\quad n_c=c(1-r),\quad s=a'+c',\quad t=a-aba.
\]
No intertwining or smallness assumption is part of this notation. When condition labels are used below, (A1) denotes $c'p=c'$, (A2) denotes $ra'=a'$, and (A3) denotes $pb=br$.
\end{notation}
\begin{lemma}[The elementary Hua identities]\label{lem:hua}
For arbitrary $a,c$ with $a+c=b^{-1}$,
\[
 t=abc=cba=c-cbc.
\]
If $a$ and $c$ are units, $t$ is a unit and $t^{-1}=a^{-1}+c^{-1}$.
\end{lemma}
\begin{proof}
Substitute $c=b^{-1}-a$ in $abc$ and $cba$. Substituting $a=b^{-1}-c$ proves the last form. When $a,c$ are units, $t=abc$ is invertible and $t^{-1}=c^{-1}b^{-1}a^{-1}=c^{-1}(a+c)a^{-1}=c^{-1}+a^{-1}$.
\end{proof}
This is the classical identity; the proof uses no commutation of $a$ and $b$. The generalization below concerns selected outer inverses, not a new proof of the elementary identity.
For comparison, if $a^-,c^-$ are inner inverses, so that
$aa^-a=a$ and $cc^-c=c$, then
\[
 t(a^-+c^-)t=t.
\]
Indeed, the two forms of $t$ give $ta^-t=tbc$ and $tc^-t=tba$,
and their sum is $tb(c+a)=t$. No intertwining is needed for this
inner-inverse identity. The issue below is instead compatibility of
selected commuting outer witnesses and inversion of the whole output.
The absorption identity $c'(a+c)a'=a'+c'$ is covered by
\cite[Theorem~8(i),(iii)]{BapatEtAl17}: for the selected outer inverses
$c',a'$, it is equivalent to $c'\R\supseteq a'\R$ and
$\R c'\subseteq\R a'$. Since $a'\R=p\R$, $\R a'=\R p$,
$c'\R=r\R$ and $\R c'=\R r$, these inclusions give $rp=p=r$.
For canonical strict ring generalized-Drazin data,
\cite[Theorem~6.1(v)]{KolihaPatricio02} gives the same equivalence by
substituting $x=-c$ and $y=a$ in
$y^{\mathrm D}-x^{\mathrm D}=x^{\mathrm D}(x-y)y^{\mathrm D}$;
the superscript here is the generalized-Drazin notation of that source.
We retain a direct proof to identify the prescribed products for the Hua
element and the ordered core--residue decomposition used below.

\pagebreak
\begin{theorem}[Exact common-core equivalence]\label{thm:common}
For the coupled data, the following are equivalent, without an intertwining or input-remainder hypothesis:
\begin{enumerate}
\item[(i)] $p=r$;
\item[(ii)] $ts=p$ and $st=r$;
\item[(iii)] $s=c'b^{-1}a'$;
\item[(iv)] $c'p=c'$ and $ra'=a'$.
\end{enumerate}
Under these conditions, write $e=p=r$ and $f=1-e$. Then $e$ commutes with $a,c,b,t,s$ and
\[
 t=T_0+\nu,\quad T_0=ete,\quad s\in e\R e,\quad sT_0=T_0s=e,
\]
\[
 \beta=fbf=(n_a+n_c)^{-1}_{f\R f},\qquad \nu=ftf=n_a\beta n_c=n_c\beta n_a=n_a-n_a\beta n_a.
\]
Moreover $s$ is the unique supported Bott--Duffin inverse of $t$ relative to $e$, and
\[
 s=e(1-e+te)^{-1},\qquad s\R=e\R,\qquad \R s=\R e.
\]
\end{theorem}
\begin{proof}
If $p=r=e$, Lemma~\ref{lem:input} gives (iv). The idempotent $e$ commutes with $a$ and $c$, hence with $b^{-1}=a+c$ and its inverse $b$.

Conversely, (iv) gives both $rp=cc'p=r$ and $rp=r(a'a)=a'a=p$, hence (i). Expanding
\[
 c'b^{-1}a'=c'(a+c)a'=c'p+ra'
\]
proves (iv)$\Rightarrow$(iii).

Conversely, if (iii) holds, multiplication on the left by $1-r$ gives $(1-r)a'=0$, and multiplication on the right by $1-p$ gives $c'(1-p)=0$. Thus (iii)$\Rightarrow$(iv).

Under (i), $b$ commutes with $e$. Using $t=abc$ and (iii),
\[
 ts=abcc'b^{-1}a'=abeb^{-1}a'=aa'=e,\qquad
 st=c'b^{-1}a'abc=c'b^{-1}ebc=c'c=e.
\]
This proves (ii).

Conversely, (ii) implies $rs=(st)s=s(ts)=sp$. Since $rc'=c'$ and $a'p=a'$, this equality becomes $ra'+c'=a'+c'p$. Multiplication by $1-r$ on the left gives $ra'=a'$, and substitution then gives $c'p=c'$. Hence (iv) follows.

All data are block diagonal along $e$; $s=ese$ and $ts=st=e$ prove the core inverse. Compressing $b(a+c)=(a+c)b=1$ by $f$ gives the inverse formula for $n_a+n_c$. The residual formulas follow by applying Lemma~\ref{lem:hua} in $f\R f$. The element $(1-e)+te$ has inverse $(1-e)+s$, proving the displayed formula for $s$. A supported Bott--Duffin inverse is the inverse of $T_0$ in $e\R e$, so uniqueness is ordinary inverse uniqueness there; this agrees with~\cite[Section 3]{Drazin12}. Finally $st=ts=e$ and $s=es=se$ give the two principal-ideal equalities.
\end{proof}
\begin{remark}[Scope of the equivalence]\label{rem:bare}
Theorem~\ref{thm:common} characterizes coincidence of the core idempotents through the prescribed product identities $ts=p$ and $st=r$; it does not characterize the bare equality $t^{\Dz}=s$. Coincidence also implies both $pb=br$ and $bp=rb$, but neither condition is imposed in that theorem. The two absorption identities already force coincidence, so they do not describe a class with distinct prescribed core idempotents.
\end{remark}

\section{Rigidity supplied by additive coupling}\label{sec:rigidity}
\begin{lemma}[Reduction under element intertwining]\label{lem:reduction}
Assume $pb=br$. Then
\[
 rp=p+n_ca',\qquad rp=r+c'n_a.
\]
Writing $z=n_cp$ and $v=rn_a$, one has
\[
 z=-n_cza',\qquad v=-c'vn_a,\qquad r-p=n_ca'-c'n_a.
\]
No remainder condition is needed.
\end{lemma}
\begin{proof}
Intertwining is equivalent to $(a+c)p=r(a+c)$. Multiplication on the right by $a'$ gives $rp=p+(1-r)ca'=p+n_ca'$. Multiplication on the left by $c'$ gives $rp=r+c'a(1-p)=r+c'n_a$. Multiplying the first relation on the left by $n_c$, using $n_cr=0$ and $pa'=a'$, gives $z=-n_cza'$. Multiplying the second on the right by $n_a$, using $pn_a=0$ and $c'r=c'$, gives the relation for $v$. Subtraction proves the last formula.
\end{proof}
\begin{theorem}[Finite-Drazin and Banach rigidity]\label{thm:rigidity}
If both input remainders are nilpotent, $pb=br$ is equivalent to $p=r$. The same equivalence holds in a complex unital Banach algebra if both input remainders are quasinilpotent.
\end{theorem}
\begin{proof}
Coincidence implies intertwining by Theorem~\ref{thm:common}.

Conversely, iteration of Lemma~\ref{lem:reduction} gives, for every $k\ge1$,
\[
 z=(-1)^k n_c^kz(a')^k,\qquad v=(-1)^k(c')^kvn_a^k.
\]
In the nilpotent case these equal zero for sufficiently large $k$. In the Banach case, if $z\ne0$, division of the norm inequality by $\|z\|$ yields $1\le\|n_c^k\|\|a'\|^k$, impossible because $\|n_c^k\|^{1/k}\to0$. The same argument with $n_a$ gives $v=0$. Since $n_ca'=za'$ and $c'n_a=c'v$, the last reduction formula gives $p=r$.
\end{proof}
The norm argument only needs vanishing root-norm limits and a submultiplicative norm, not completeness. Outside a Banach algebra we make no spectrum-based identification of this hypothesis.
\begin{remark}[The classical intertwiner mechanism]\label{rem:intertwiner}
Put $d=p-r$. Expanding $(a+c)p=r(a+c)$ gives $cd=-da$. Drazin's classical result transfers intertwiners to finite Drazin inverses~\cite[Theorem 2.2]{Drazin13}. Applying it to $c,-a$, and $d$ gives $c'd=-da'$ in the finite-Drazin case. Since $c'd+da'=c'p-c'+a'-ra'$, multiplying on the left by $1-r$ and on the right by $1-p$ yields the absorptions in Theorem~\ref{thm:common}. This supplies another proof of coincidence. The coupling-specific content is the relation $cd=-da$, not the inverse-intertwining principle itself. The direct iteration above also displays the one-sided obstruction when the finite or Banach hypotheses are absent.

The rigidity conclusion needs only $(a+c)p=r(a+c)$: for arbitrary
$a,c$ with commuting outer inverses, this identity implies $p=r$ if
both input remainders are nilpotent, or if they are quasinilpotent in a
complex unital Banach algebra, even when $a+c$ is not a unit.
The reduction equations and their finite iterations or root-norm
estimates use that identity alone. This observation concerns input
coincidence. The Hua element and residual parallel sum in
Theorem~\ref{thm:common} retain the unit hypothesis of
Notation~\ref{not:coupled}.

Without additive coupling, unit intertwining alone does not force coincidence.
In $M_2(\C)$, take $a=a'=p=E_{11}$, $c=c'=r=E_{22}$, and
$b=E_{12}+E_{21}$. The input remainders vanish and $pb=br$, but $p\ne r$;
here $a+c=1\ne b^{-1}$.
\end{remark}
\begin{theorem}[Separated one-sided power filtrations]\label{thm:separated}
Assume $pb=br$ and
\[
 \bigcap_{k\ge1}n_c^k\R=\{0\},\qquad \bigcap_{k\ge1}\R n_a^k=\{0\}.
\]
Then $p=r$. No norm, commutativity of the ring, or a priori quasinilpotence is required.

In particular, suppose $\R$ has a decreasing filtration by additive subgroups
\[
 \R=F^0\supseteq F^1\supseteq\cdots,\qquad
 F^jF^k\subseteq F^{j+k}\quad(j,k\ge0),\qquad
 \bigcap_{k\ge1}F^k=\{0\}.
\]
If $n_a,n_c\in F^1$, then $p=r$.
\end{theorem}
\begin{proof}
The iterates in the preceding proof put $z$ in every $n_c^k\R$ and $v$ in every $\R n_a^k$. Separation forces both to vanish; apply the reduction formula.

For the filtration consequence, induction gives $n_c^k,n_a^k\in F^k$.
Multiplicativity with $F^0=\R$ gives $n_c^k\R\subseteq F^k$ and
$\R n_a^k\subseteq F^k$; hence the two one-sided intersections vanish.
\end{proof}

The multiplicative-filtration special case is necessarily classical.
Indeed, after coincidence put $f=1-e$ and $u=n_a+n_c$. By
Theorem~\ref{thm:common}, $u$ is invertible in $f\R f$.
Since $u\in F^1$ and every $F^k$ is a two-sided ideal,
\[
 f=u^k(u^{-1}_{f\R f})^k\in F^k\qquad(k\ge1).
\]
Thus $f=0$. This restriction applies to the filtration special case,
not to the general one-sided separation criterion.

These are explicit one-sided hypotheses. The terminology of commutative adic separation provides a comparison with Atiyah and Macdonald~\cite[Theorem 10.17 and Corollaries 10.19--10.20]{AtiyahMacdonald69}. The intersection theorem assumes a commutative Noetherian ring and a finitely generated module. Corollary 10.19 gives vanishing when the ideal lies in the Jacobson radical; Corollary 10.20 applies to the maximal ideal of a Noetherian local ring. Here the one-sided separation hypotheses are assumed directly, and the ring need not be commutative or Noetherian.

\begin{proposition}[Coincidence modulo the radical]\label{prop:radical}
Assume $pb=br$. If some power of each input remainder lies in $J(\R)$, then
\[
 z=n_cp\in J(\R),\qquad v=rn_a\in J(\R),\qquad p-r\in J(\R).
\]
Moreover the element
\[
 w=rp+(1-r)(1-p)
\]
belongs to $1+J(\R)$ and satisfies $r=wpw^{-1}$.
\end{proposition}
\begin{proof}
Choose $k\ge1$ so that $n_c^k,n_a^k\in J(\R)$. The iterates in the proof of
Theorem~\ref{thm:rigidity} give $z,v\in J(\R)$. The identity
$r-p=za'-c'v$ then gives $p-r\in J(\R)$. Equivalently, the complete equations
pass to $\R/J(\R)$ with nilpotent input remainders and a unit image of $b$,
and Theorem~\ref{thm:rigidity} gives $\bar p=\bar r$ there.
Consequently $\bar w=\bar p^2+(1-\bar p)^2=1$, so
$w\in1+J(\R)\subseteq\R^\times$. Finally,
\[
 wp=rp=rw,
\]
which gives the asserted conjugacy. These conclusions do not imply $p=r$ in
$\R$.
\end{proof}

The radical conclusion strengthens to actual coincidence when
$\bigcap_{j\ge1}J(\R)^j=0$. If powers of both input remainders
belong to $J(\R)$, sufficiently long iterates of Lemma~\ref{lem:reduction}
place $z$ and $v$ in every radical power, forcing both to vanish.

The next theorem supplies a structural hypothesis under which the
one-sided intersections in Theorem~\ref{thm:separated} vanish.

\begin{theorem}[Separation over a Noetherian central subring]
\label{thm:central-separation}
Let $\R$ be a unital ring that is finitely generated as a module over
a commutative Noetherian subring $Z$ of its center, with the same
identity. For every Harte-quasinilpotent $q\in\R$,
\[
 \bigcap_{k\ge1}q^k\R=\{0\},\qquad
 \bigcap_{k\ge1}\R q^k=\{0\}.
\]
Consequently, for the coupled data of Notation~\ref{not:coupled} with
Harte-quasinilpotent input remainders,
\[
 pb=br\quad\Longleftrightarrow\quad p=r.
\]
\end{theorem}
\begin{proof}
The zero ring is immediate, so assume $\R\ne0$ and put $S=Z[q]$.
This is a commutative subring of $\R$. It is finite over $Z$, being a
$Z$-submodule of the finite module $\R$ over the Noetherian ring $Z$;
hence $S$ is Noetherian. The same finite set of $Z$-generators makes
$\R$ a finite left and a finite right $S$-module.

The subring $S$ is inverse-closed in $\R$. Indeed, if
$u\in S\cap\R^\times$, then $u\R=\R$. The determinant form of
Nakayama's lemma, applied to the finite left $S$-module $\R$, gives
$v\in S$ such that $(1-uv)\R=0$
\cite[Corollary~2.5, p.~21]{AtiyahMacdonald69}. Evaluation at $1$
gives $uv=1$; commutativity of $S$ makes $u$ a unit of $S$.
For every $s\in S$, Harte quasinilpotence makes $1-qs$ a unit of
$\R$, and hence of $S$. Thus $q\in J(S)$.

Krull intersection for the ideal $(q)\subseteq J(S)$ on the two finite
$S$-module structures now gives both asserted intersections
\cite[Theorem~10.17 and Corollary~10.19, pp.~110--111]{AtiyahMacdonald69}.
Apply Theorem~\ref{thm:separated} to the input remainders; the reverse
implication follows from Theorem~\ref{thm:common}.
\end{proof}
The separation conclusion does not imply nilpotence: $z$ is
Harte-quasinilpotent and nonnilpotent in $D=\C[z]_{(z)}$, while
$\bigcap_{j\ge1}z^jD=0$.

The finite-fiber observation needs no Noetherianity. If $\R$ is finite
over a commutative central subring $Z$ with the same identity and
$q\in\qnil(\R)$, then $\bar q$ is nilpotent in every
$\R/\mathfrak m\R$, where $\mathfrak m$ is maximal in $Z$.
For a nonzero fiber, otherwise its minimal polynomial is
$X^r h(X)$ with $h(0)\ne0$ and $h$ nonconstant, so $h(\bar q)$
is not a unit. Writing $h(X)/h(0)=1+Xg(X)$ and lifting the
coefficients of $g$ to $Z$ contradicts the Harte unit test for the
resulting polynomial in $q$. The zero fiber is immediate. This
argument uses lifted polynomial tests, not general preservation of
Harte quasinilpotence under quotients.

\begin{theorem}[Intertwining over a finite central subring]
\label{thm:central-transfer}
Let $\R$ be finite as a module over a commutative central subring $Z$
with the same identity. Let $x,y\in\R$ have commuting outer inverses
$x',y'$, put $p_x=xx'$, $p_y=yy'$, and assume
\[
 n_x=x(1-p_x)\in\qnil(\R),\qquad
 n_y=y(1-p_y)\in\qnil(\R).
\]
For every $h\in\R$,
\[
 xh=hy\quad\Longrightarrow\quad
 p_xh=hp_y\quad\hbox{and}\quad x'h=hy'.
\]
Consequently, for the coupled data of Notation~\ref{not:coupled}
with Harte-quasinilpotent input remainders, $pb=br$ is equivalent
to $p=r$. No Noetherian hypothesis on $Z$ is required.
\end{theorem}
\begin{proof}
Every $v\in\R$ is integral over $Z$: apply the determinant trick to
left multiplication by $v$ on the finite $Z$-module $\R$, then
evaluate its monic annihilating polynomial at $1$
\cite[Proposition~2.4, p.~21]{AtiyahMacdonald69}.

We first observe that, for $n\in\qnil(\R)$ and $v$ integral over $Z$,
either $u=nuv$ or $u=vun$ forces $u=0$. Choose
\[
 v^m+c_{m-1}v^{m-1}+\cdots+c_0=0,\qquad c_i\in Z.
\]
In the first case $u=n^juv^j$ for every $j\ge0$, whence
\[
 u=n^muv^m=-\sum_{i=0}^{m-1}c_i n^muv^i
             =-\sum_{i=0}^{m-1}c_i n^{m-i}u.
\]
Thus $(1+c_{m-1}n+\cdots+c_0n^m)u=0$.
This factor is a unit by Harte quasinilpotence, since it has the form
$1+ng(n)$ with central polynomial coefficients. The second case
puts the same unit on the right. No commutation of $n$ with $v$
is used.

Set $X=(1-p_x)hp_y$, $Y=p_xh(1-p_y)$, and $Z_0=p_xhp_y$.
Compression of $xh=hy$ gives
\[
 n_xX=X(yp_y),\qquad (xp_x)Y=Yn_y,\qquad
 (xp_x)Z_0=Z_0(yp_y).
\]
Multiplication by the corner inverses yields $X=n_xXy'$ and
$Y=x'Yn_y$. The preceding observation gives $X=Y=0$, so
$p_xh=Z_0=hp_y$. Multiplying the third compressed equality on
the left by $x'$ and on the right by $y'$ gives
$x'Z_0=Z_0y'$, and therefore $x'h=hy'$.

For the coupled data put $d=p-r$. Intertwining gives $cd=-da$.
Apply the transfer just proved to $(x,y,h)=(c,-a,d)$ with
witnesses $(c',-a')$. The remainder $-n_a$ is Harte-quasinilpotent,
and support transfer gives $r(p-r)=(p-r)p$, or $2rp=p+r$.
Multiplication on the left by $r$ and on the right by $p$ gives,
respectively, $rp=r$ and $rp=p$, without division by $2$.
Thus $p=r$; the converse is Theorem~\ref{thm:common}.
\end{proof}
The same corner argument transfers intertwiners when
$\bigcap_k n_x^k\R=0$ and $\bigcap_k\R n_y^k=0$:
iteration puts $X$ and $Y$ in the respective intersections.
It also gives the finite-Drazin and Banach root-norm versions.
Input coincidence and output bicommutation are applications of this
mechanism, but their failures in Section~\ref{sec:ring} are logically
different: the second occurs after the input cores already coincide.

\emph{Power separation is stronger.}
The Noetherian hypothesis in Theorem~\ref{thm:central-separation}
cannot simply be omitted. Let $D=\C[z]_{(z)}$, $K=\C(z)$, and
$S=D\ltimes K$, with multiplication
\[
 (d,k)(d',k')=(dd',dk'+d'k).
\]
This commutative ring is finite over itself, but the ideals
$0\ltimes z^{-j}D$ form a strictly ascending chain. A pair $(d,k)$
is a unit exactly when $d$ is a unit of $D$, so
$J(S)=zD\ltimes K$. Thus $q=(z,0)$ is Harte-quasinilpotent, whereas
\[
 q^jS=z^jD\ltimes K,\qquad
 \bigcap_{j\ge1}q^jS=0\ltimes K\ne0.
\]
Theorem~\ref{thm:central-transfer} eliminates a recurrence with one
fixed integral coefficient; it does not assert universal power
separation.

\section{The general-ring boundaries}\label{sec:ring}
The preceding positive results do not imply coincidence for Harte-quasinilpotent input remainders in arbitrary rings. Nor does coincidence alone ensure strict generalized-Drazin invertibility of the output. The following two constructions separate these questions. All inverses in the displayed localized rings are taken in those rings, not just in a larger matrix algebra over their fraction fields.
\begin{example}[Failure of general-ring coincidence]\label{ex:ring}
Let $D=\C[z]_{(z)}$, the rational functions regular at zero, let $K=\C(z)$, and set
\[
 R=\left\{\begin{pmatrix}H&v\\0&k\end{pmatrix}:H\in M_2(D),\ v\in K^2,\ k\in K\right\}.
\]
The choice of the two coefficient rings is essential: $z$ is not a unit
of $D$, whereas the off-diagonal module $K^2$ permits division by $z$.
This allows a nonzero obstruction to persist under the one-sided
recurrence of Lemma~\ref{lem:reduction}. Take
\[
 a=\begin{pmatrix}0&0&-z\\1&0&-1\\0&0&-z\end{pmatrix},\qquad
 c=\begin{pmatrix}z&1&0\\0&z&0\\0&0&1\end{pmatrix},\qquad
 a'=-z^{-1}(E_{13}+E_{33}),\quad c'=E_{33}.
\]
Then
\[
 B:=a+c=\begin{pmatrix}z&1&-z\\1&z&-1\\0&0&1-z\end{pmatrix},\qquad
 b=B^{-1}=\begin{pmatrix}
 \dfrac{z}{z^2-1}&\dfrac{-1}{z^2-1}&\dfrac{1}{1-z}\\
 \dfrac{-1}{z^2-1}&\dfrac{z}{z^2-1}&0\\
 0&0&\dfrac{1}{1-z}
 \end{pmatrix}.
\]
The top-left determinant $z^2-1$ is a unit of $D$, and $1-z\ne0$ in $K$; hence $B,b\in R^\times$. Multiplication gives
\[
 a'aa'=a',\quad aa'=a'a=p=E_{13}+E_{33},\qquad
 c'cc'=c',\quad cc'=c'c=r=E_{33}.
\]
In particular $p-r=E_{13}\ne0$. The input remainders are
\[
 n_a=\begin{pmatrix}0&0&0\\1&0&-1\\0&0&0\end{pmatrix},\qquad
 n_c=\begin{pmatrix}z&1&0\\0&z&0\\0&0&0\end{pmatrix}.
\]
Here $n_a^2=0$ and $n_c$ is not nilpotent. To verify its smallness inside $R$, note that
\[
 J(R)=\left\{\begin{pmatrix}H&v\\0&0\end{pmatrix}:H\in M_2(zD),\ v\in K^2\right\}.
\]
Indeed, the displayed set is an ideal $I$ and $1-yx$ is a unit of $R$ for every $x\in I,y\in R$: its top-left block is congruent to $I_2$ modulo $zD$ and its lower diagonal entry is $1$. Thus $I\subseteq J(R)$; the quotient $M_2(\C)\times K$ has zero radical, proving equality. Consequently $n_c^2\in J(R)$. If $q^2\in J(R)$ and $x$ commutes with $q$, then $1-(qx)^2$ is a unit. Therefore, its commuting factors $1+qx$ and $1-qx$ are units, proving that $q$ is Harte-quasinilpotent. Both input remainders have that property here.

The primed elements satisfy the strict double-commutant convention. For $a$, put $A=E_{21}\in M_2(D)$ and $S=I_3+E_{13}$. Then
\[
 a=S\diag(A,-z)S^{-1},\qquad a'=S\diag(0_2,-z^{-1})S^{-1}.
\]
An element $\left(\begin{smallmatrix}H&v\\0&k\end{smallmatrix}\right)$ commuting with $\diag(A,-z)$ satisfies $(A+zI_2)v=0$. The coefficient is invertible over $K$, so $v=0$, and the displayed inverse commutes with every such element. Conjugating back gives $a'\in\comm^2(a)$. For $c=\diag(zI_2+E_{12},1)$ the corresponding equation is $((z-1)I_2+E_{12})v=0$, again forcing $v=0$. Thus $c'\in\comm^2(c)$. The primed elements are strict ring generalized Drazin inverses, and also satisfy the pseudo-Drazin conditions of Wang and Chen~\cite[Definition 3.1, p.~1338]{WangChen12}, since both input-remainder squares lie in $J(R)$.

Finally,
\[
 pb=br=\frac{1}{1-z}p,\qquad Bp=rB=(1-z)r.
\]
Even operator intertwining holds. For $X=\left(\begin{smallmatrix}H&v\\0&k\end{smallmatrix}\right)\in R$, one has $pXp=kp$, $rXr=kr$, and $brb=p/(1-z)^2$. The lower diagonal entry of $bXb$ is $k/(1-z)^2$, whence
\[
 p(bXb)p=\frac{k}{(1-z)^2}p=b(rXr)b.
\]
The obstruction in Lemma~\ref{lem:reduction} is $n_cp=zE_{13}\ne0$, whereas $rn_a=0$. Its required cross-commutation fails: $n_ca'=-E_{13}$ but $a'n_c=0$. Thus strict generalized-Drazin input hypotheses do not force coincidence in an arbitrary ring, and the radical conclusion of Proposition~\ref{prop:radical} cannot be strengthened to equality. The nonnilpotent remainder $n_c$ does not contradict finite-Drazin rigidity, and the example does not satisfy the absorption identities equivalent to coincidence in Theorem~\ref{thm:common}.
\end{example}
\begin{example}[Failure of strict output inversion with a common core]\label{rem:ring-output}
We construct finite-Drazin input data with coincident nonzero proper
cores and a Harte-quasinilpotent residue, but no strict generalized
Drazin inverse for the output. The obstruction is a nonzero
off-diagonal element commuting with the output but not with its
supported candidate.

\emph{The strict inverse criterion.}
In a general ring, membership of the candidate in the double commutant is a further condition, and it is exactly what separates the Bott--Duffin inverse relations from the generalized Drazin inverse of~\cite{KolihaPatricio02}. Suppose $st=ts=e=e^2$ and $sts=s$, and write $\nu=t(1-e)$, as produced by Theorem~\ref{thm:common}. Then $s$ is the generalized Drazin inverse of $t$ in the double-commutant sense if and only if $\nu$ is Harte quasinilpotent and $e\in\comm^2(t)$, equivalently $s\in\comm^2(t)$. The inverse-axiom remainder is $t-t^2s=t(1-e)=\nu$. If $s\in\comm^2(t)$, then $e=ts\in\comm^2(t)$.

Conversely, if $e\in\comm^2(t)$, then $1-e+te$ is invertible with inverse $1-e+s$, so any element commuting with $t$ and with $e$ commutes with $s=e(1-e+te)^{-1}$. The extra condition is not automatic, even for complete Hua data.

\emph{The coupled construction.} In the ring $\left(\begin{smallmatrix}M_3(D)&K^3\\0&K\end{smallmatrix}\right)$ with $D=\C[z]_{(z)}$ and $K=\C(z)$, take
\[
 N=E_{12}+E_{23},\qquad C_z=-zE_{11}+E_{12}-z^2E_{21}+zE_{22}+E_{31}.
\]

Put
\[
 B_z:=\begin{pmatrix}
 0&0&1\\[1pt]
 \tfrac12&0&\tfrac z2\\[1pt]
 -\tfrac z2&1&\tfrac{z^2}2
 \end{pmatrix}\in M_3(D),\qquad
 N+C_z=\begin{pmatrix}-z&2&0\\-z^2&z&1\\1&0&0\end{pmatrix}.
\]
Direct multiplication gives
\[
 (N+C_z)B_z=B_z(N+C_z)=I_3,\qquad
 NB_zC_z=F_z:=\tfrac12E_{12}+\tfrac z2E_{22}.
\]
Now take
\[
 a=\diag(N,z),\quad c=\diag(C_z,z),\quad
 b=\diag(B_z,\tfrac1{2z}),\quad a'=c'=\diag(0_3,z^{-1}).
\]
These are complete Hua data satisfying (A1)--(A3), with common idempotent $e=\diag(0_3,1)$, $0\ne e\ne1$, Hua element $t=\diag(F_z,z/2)$, candidate $s=\diag(0_3,2/z)$, and residue $\nu=\diag(F_z,0)$. Both input indices are exactly three: $N^2=E_{13}\ne0$, $N^3=0$, and
\[
 C_z=(I_3+zE_{21})C_0(I_3-zE_{21}),\qquad
 C_0=E_{12}+E_{31},\qquad C_z^2=E_{32}-zE_{31}\ne0,\quad C_z^3=0.
\]
For $H\in M_3(D)$ commuting with $F_z$, reduction modulo $zD$ gives $\overline F_z^{\,2}=0$ and $(\overline F_z\,\overline H)^2=0$. Hence $\det(I_3+F_zH)$ reduces to $1$ modulo $zD$ and is a unit of the local ring $D$. The adjugate formula gives an inverse in $M_3(D)$, proving Harte quasinilpotence of $F_z$ there. It is not nilpotent, since $F_z^k=(z/2)^{k-1}F_z\ne0$ for every $k\ge1$.

For the whole triangular ring, let $Z=\left(\begin{smallmatrix}H&v\\0&k\end{smallmatrix}\right)$ commute with $\nu=\diag(F_z,0)$. Its upper-left block satisfies $HF_z=F_zH$. The diagonal blocks of $1+\nu Z$ are $I_3+F_zH$ and $1$, both invertible, so the triangular inverse belongs to the specified ring. Hence $\nu$ is Harte-quasinilpotent there.

\emph{The off-diagonal obstruction.} Let
\[
 w=(1,z,0)^{\mathsf T},\qquad
 X=\begin{pmatrix}0_3&w\\0&0\end{pmatrix}\in\R.
\]
Since $F_zw=(z/2)w$, one has $Xt=tX$. However, $sX=0$ and $Xs=(2/z)X\ne0$, so $s\notin\comm^2(t)$. To exclude every other strict generalized Drazin inverse, suppose that $W\in\comm^2(t)$, $WtW=W$, and $t-t^2W$ is Harte-quasinilpotent in $\R$. The idempotent $e$ commutes with $t$, so $W$ commutes with $e$ and has the form $W=\diag(Y,\eta)$ with $Y\in M_3(D)$ and $\eta\in K$. Commutation with $t$ and the outer equation give
\[
 F_zY=YF_z,\qquad YF_zY=Y.
\]
Thus $P:=F_zY$ is idempotent. Harte quasinilpotence of $F_z$ and commutation of $Y$ with $F_z$ make $I_3-P$ a unit in $M_3(D)$. Since $P(I_3-P)=0$, it follows that $P=0$, and then $Y=YF_zY=YP=0$. The actual remainder is
\[
 t-t^2W=\diag\bigl(F_z-F_z^2Y,\ z/2-(z^2/4)\eta\bigr).
\]
Its compression to $e\R e\cong K$ is Harte-quasinilpotent. The only Harte-quasinilpotent element of a field is zero, so $z/2-(z^2/4)\eta=0$ and $\eta=2/z$. Consequently $W=s$, contradicting $s\notin\comm^2(t)$. Thus $t$ has no generalized Drazin inverse in the strict ring sense.

\emph{Commuting outer inverses without bicommutation.}
For every $\lambda\in\C$, the element $s+\lambda X$ is nevertheless a commuting outer inverse of $t$, and $(s+\lambda X)t=e+\lambda(z/2)X$. Its remainder is Harte-quasinilpotent as well: with $\mu=z/2$, the unit $U_\lambda=1+\lambda\mu X$ commutes with $t$ and satisfies
\[
 U_\lambda^{-1}=1-\lambda\mu X,\qquad
 U_\lambda sU_\lambda^{-1}=s+\lambda X,\qquad
 t-t^2(s+\lambda X)=U_\lambda\nu U_\lambda^{-1}.
\]
The prescribed idempotent changes when $\lambda\ne0$. Therefore, these inverses do not contradict uniqueness relative to the fixed $e$, or uniqueness of $s=a'+c'$ for the given input tuple. Example~\ref{ex:product} satisfies the double-commutant condition because its candidate is central in the product ring. Coincidence and output double-commutant membership are separate conditions.
\end{example}
\begin{remark}[Two logically separate obstructions]\label{rem:two}
In the first construction the input witnesses satisfy the strict double-commutant conditions, but the core idempotents differ. In the second, the common core is already present and the difficulty is an output double commutant. A formula valid for a supported outer inverse does not by itself overcome either obstruction. Therefore, the Banach block formula below has a Banach hypothesis, and the finite-Drazin block formula has its own algebraic proof.
\end{remark}
\begin{proposition}[Where the noncommutativity must live]\label{prop:cross}
Assume Notation~\ref{not:coupled} and (A3).
\begin{enumerate}
\item[(i)] If $n_c\in\qnil(\R)$ and $[n_c,a']=0$, then $z=n_cp=0$.
\item[(ii)] If $n_a\in\qnil(\R)$ and $[n_a,c']=0$, then $v=rn_a=0$.
\item[(iii)] Consequently, if both input remainders are quasinilpotent and both cross-commutations hold, then $p=r$.
\end{enumerate}
\end{proposition}
\begin{proof}
(i) The absorption $pa'=a'$ (Lemma~\ref{lem:input}) gives $za'=n_cpa'=n_ca'$, which equals $a'n_c$ by hypothesis. Therefore, the fixed-point equation $z=-n_cza'$ in Lemma~\ref{lem:reduction} evaluates in closed form:
\[
 z=-n_c(za')=-n_c(a'n_c)=-n_c^2a',
\]
while $zn_ca'=z(a'n_c)=(za')n_c=(n_ca')n_c=n_c^2a'$. Adding,
\[
 z(1+n_ca')=-n_c^2a'+n_c^2a'=0.
\]
Since $a'\in\comm(n_c)$ and $n_c$ is quasinilpotent, $1+n_ca'$ is a unit directly by Harte's definition, and right cancellation gives $z=0$.

(ii) This is the mirror image of part (i): $c'=c'r$ gives $c'v=c'rn_a=c'n_a=n_ac'$, hence $v=-(c'v)n_a=-n_ac'n_a=-n_a^2c'$ and $n_ac'v=n_a^2c'$. Thus $(1+n_ac')v=0$; Harte quasinilpotence makes $1+n_ac'$ a unit, so left cancellation gives $v=0$.

(iii) Parts (i) and (ii) give $z=v=0$. Since $n_ca'=za'$ and $c'n_a=c'v$, the final identity in Lemma~\ref{lem:reduction} gives $r-p=0$.
\end{proof}
\begin{remark}[Mirrored intertwining]\label{rem:mirror}
Passing to the opposite ring exchanges left and right in the preceding arguments. Therefore, the corresponding conclusions hold for $bp=rb$ under the same two-sided smallness hypotheses.
\end{remark}

The off-diagonal recurrence also controls the independent output
obstruction in Example~\ref{rem:ring-output}.

\begin{proposition}[Separation and strict output inversion]
\label{prop:separated-output}
Let $e\in\R$ be an idempotent, put $f=1-e$, and suppose
\[
 t=T_0+\nu,\qquad T_0\in(e\R e)^\times,\qquad
 \nu\in f\R f,\qquad s=(T_0)^{-1}_{e\R e}.
\]
If
\[
 \bigcap_{k\ge1}\nu^k\R=\{0\},\qquad
 \bigcap_{k\ge1}\R\nu^k=\{0\},
\]
then $e,s\in\comm^2(t)$.

For common-core coupled data, if $\R$ is finite as a module over a
commutative central subring with the same identity, then $s$ is the
strict ring generalized Drazin inverse of $t$ if and only if $\nu$
is Harte-quasinilpotent in $f\R f$. This conclusion does not require
Noetherianity of the central subring.
\end{proposition}
\begin{proof}
For $x\in\comm(t)$ set $u=exf$ and $v=fxe$. Compression of $tx=xt$
gives $T_0u=u\nu$ and $\nu v=vT_0$. For every $k\ge1$,
\[
 u=s^ku\nu^k\in\R\nu^k,\qquad
 v=\nu^kvs^k\in\nu^k\R.
\]
Separation gives $u=v=0$. Thus every $x\in\comm(t)$ commutes with
$e$, proving $e\in\comm^2(t)$. It then commutes with
$1-e+te$ and its inverse, and hence with
$s=e(1-e+te)^{-1}$. This proves $s\in\comm^2(t)$.

For the last assertion, a Harte-quasinilpotent element
$\nu\in f\R f$, extended by zero on $e$, is Harte-quasinilpotent
in $\R$. Indeed, if $x\nu=\nu x$, the off-diagonal products
$exf\nu$ and $\nu fxe$ vanish, while $y=fxf$ commutes with
$\nu$. Hence $1+\nu x=e+(f+\nu y)$ is invertible, with inverse
$e+(f+\nu y)^{-1}_{f\R f}$. Since $sts=s$, $st=ts=e$, and
$t-t^2s=\nu$, apply Theorem~\ref{thm:central-transfer} to the pair
$(t,t)$ with both witnesses equal to $s$. Every element commuting
with $t$ then commutes with $s$. Thus $s$ is the strict inverse. Conversely, a strict inverse $s$ has
Harte-quasinilpotent remainder $\nu=t-t^2s$ in $\R$; compression
to $f\R f$ preserves that property by Lemma~\ref{lem:qnil}(b).
The zero-corner cases use the convention in Section~\ref{sec:common}.
\end{proof}
The residual quasinilpotence hypothesis is separate from the input
hypotheses. The finite-module condition removes the double-commutant
obstruction; it does not force the residual parallel sum to be
Harte-quasinilpotent.

\begin{corollary}[Bicommutation without Hua coupling]
\label{cor:central-bicommutation}
Let $\R$ be finite as a module over a commutative central subring
with the same identity. If $x,y\in\R$ satisfy
\[
 xy=yx,\qquad yxy=y,\qquad x-x^2y\in\qnil(\R),
\]
then $y\in\comm^2(x)$, so $y$ is the strict ring generalized Drazin
inverse of $x$. No Noetherian hypothesis is required.
\end{corollary}
\begin{proof}
For every $h\in\comm(x)$, apply Theorem~\ref{thm:central-transfer}
to the pair $(x,x)$ with both selected witnesses equal to $y$.
It gives $yh=hy$.
\end{proof}
The corollary assumes the commuting outer witness and its small
remainder; it does not assert inverse existence for every element.

For comparison, the rings in Examples~\ref{ex:ring}
and~\ref{rem:ring-output} have the form
\[
 R_d=\begin{pmatrix}M_d(D)&K^d\\0&K\end{pmatrix},\qquad
 D=\C[z]_{(z)},\quad K=\C(z),\quad d=2,3.
\]
Their local diagonal ring and fraction-field connecting module use
the mechanism of Herstein's one-sided Noetherian example, recalled by
Brown and Lenagan~\cite[p.~7]{BrownLenagan82}. The constructions above
add the specified coupled inverse data to that mechanism.
Their center is $\{\diag(d_0I_d,d_0):d_0\in D\}$: commutation
with the diagonal corners eliminates the off-diagonal block, and
commutation with every connecting column identifies the two scalar
diagonal entries. The module $K^d$ is not finite over $D$, since a
finite set of generators has bounded pole order at zero. Thus $R_d$
is not finite over its center, or over any smaller unital central
subring. The radical calculation in Example~\ref{ex:ring} gives
\[
 J(R_d)^j=\begin{pmatrix}M_d(z^jD)&K^d\\0&0\end{pmatrix}
 \quad(j\ge1),\qquad
 \bigcap_{j\ge1}J(R_d)^j
 =\begin{pmatrix}0&K^d\\0&0\end{pmatrix},
\]
because $zK^d=K^d$ and $\bigcap_jz^jD=0$. The connecting module
therefore survives every radical power.

These rings are right Noetherian. With $\epsilon=\diag(I_d,0)$,
the right module $\epsilon R_d$ has submodule
$\epsilon R_d(1-\epsilon)\cong K^d$ and quotient $M_d(D)$;
both are Noetherian as right $R_d$-modules, as is
$(1-\epsilon)R_d\cong K$. In contrast, the upper-right left ideals
with columns in $z^{-j}D^d$ form a strictly ascending chain.
Passing to the opposite ring and interchanging $a,c$ gives
left-Noetherian counterexamples. Thus either one-sided Noetherian
hypothesis alone is insufficient; no conclusion about two-sided
Noetherianity alone is asserted.

\section{Core and residual inverse criteria}\label{sec:inverses}
In this section assume $p=r=e$ and put $f=1-e$, $B=f\R f$. Theorem~\ref{thm:common} supplies $t=T_0+\nu$, with $T_0$ invertible in $e\R e$ and inverse $s$. It does not impose any spectral or nilpotence condition on $\nu$.
\begin{lemma}[Block inversion]\label{lem:blocks}
If $x$ commutes with an idempotent $e$, then $x$ is invertible in $\R$ if and only if its two diagonal corners are invertible in their respective corner rings. In that case its inverse is the sum of their inverses. For a complex unital Banach algebra,
\[
 \sigma_\A(x)=\sigma_{e\A e}(exe)\cup\sigma_{f\A f}(fxf).
\]
\end{lemma}
\begin{proof}
If $x$ is invertible, its inverse commutes with $e$ and its compressions are the required corner inverses.

Conversely, the sum of the two corner inverses multiplies $x$ to $e+f=1$ on both sides. Apply this equivalence to $\lambda1-x$ for the spectral assertion. Empty corner spectra use the convention fixed above.
\end{proof}
\begin{theorem}[Output existence and success of the candidate]\label{thm:output}
With a common core idempotent as above:
\begin{enumerate}
\item[(a)] $t$ is group invertible if and only if $\nu$ is group invertible in $B$, and then $t^\#=s+\nu_B^\#$.
\item[(b)] $t$ is finite-Drazin invertible if and only if $\nu$ is finite-Drazin invertible in $B$, and then $t^{\Dz}=s+\nu_B^{\Dz}$.
\item[(c)] $s=t^\#$ if and only if $\nu=0$, and $s=t^{\Dz}$ if and only if $\nu$ is nilpotent. When $s=t^{\Dz}$ and $e\ne1$, $\ind(t)$ equals the nilpotency index of $\nu$. If $e=1$, then $t$ is a unit and $\ind(t)=0$.
\item[(d)] The element $s$ is always group invertible, with $s^\#=T_0=t^2s$. Thus $s^{\Dz}=t$ if and only if $\nu=0$.
\end{enumerate}
These statements include the zero-core and zero-residual-corner cases.
\end{theorem}

\begin{proof}
By Theorem~\ref{thm:common},
\[
 t=T_0+\nu,\qquad T_0,s\in e\R e,\qquad \nu\in B=f\R f,
 \qquad T_0s=sT_0=e.
\]
Products between the two corners vanish. In particular,
$ts=st=e$, $sts=s$, and $t-t^2s=\nu$.

For clarity, if $y$ is a commuting outer inverse of an element $x$
in a unital ring, then $xy$ is an idempotent commuting with $x$, and
\[
 (x-x^2y)^m=x^m(1-xy)=x^m-x^{m+1}y\qquad(m\ge1).
\]
Thus the nilpotent-remainder condition is equivalent to the
finite-Drazin power equation for some positive exponent.

\emph{(a) and (b).}
If $v$ is the group or finite Drazin inverse of $\nu$ in $B$,
set $w=s+v$. Orthogonality of the corners gives
\[
 tw=wt=e+\nu v,\qquad
 wtw=sT_0s+v\nu v=s+v=w,
 \qquad t-t^2w=\nu-\nu^2v.
\]
The last expression is zero in the group case and nilpotent in the
finite-Drazin case. Positive powers of an element of $B$ are unchanged
when computed in $\R$, so $w$ is the corresponding inverse of $t$.

Conversely, let $w$ be the group or finite Drazin inverse of $t$.
By~\cite[Theorem~1]{Drazin58}, $w\in\comm^2(t)$.
Since $e\in\comm(t)$, we have $we=ew$; hence
\[
 w=u+v,\qquad u=ewe\in e\R e,\qquad v=fwf\in B.
\]
Compressing $tw=wt$ and $wtw=w$ gives
\[
 T_0u=uT_0,\qquad uT_0u=u,\qquad
 \nu v=v\nu,\qquad v\nu v=v.
\]
In the group case, compressing $twt=t$ gives
$T_0uT_0=T_0$ and $\nu v\nu=\nu$.
Multiplying the first equality by $s$ on both sides gives $u=s$;
the remaining equations say that $v=\nu_B^\#$.

In the finite-Drazin case, there is a $k\ge1$ such that
$t^{k+1}w=t^k$. Compression gives
\[
 T_0^{k+1}u=T_0^k,\qquad \nu^{k+1}v=\nu^k.
\]
Multiplying the first equality on the left by $s^{k+1}$ gives $u=s$.
The equations for $v$ show that $v=\nu_B^{\Dz}$.
Thus $w=s+v$ in each case, proving both converse implications and
the inverse formulas.

\emph{(c).}
The identities $sts=s$, $st=ts=e$, and $t-t^2s=\nu$ give
$s=t^\#$ exactly when $\nu=0$, and $s=t^{\Dz}$ exactly when
$\nu$ is nilpotent.
Assume $s=t^{\Dz}$. For every $k\ge1$, orthogonality gives
$t^kf=\nu^k$, and therefore,
\[
 t^{k+1}s=t^k(ts)=t^ke=t^k-\nu^k.
\]
Consequently the least admissible positive exponent is the
nilpotency index of $\nu$. The exponent-zero equation is
$ts=1$, which is equivalent to $e=1$.
Hence, if $e\ne1$, $\ind(t)$ is the nilpotency index of $\nu$,
including index one when $\nu=0$.
If $e=1$, then $\nu=0$ and $ts=st=1$, so $t$ is a unit and
$\ind(t)=0$.

\emph{(d).}
The identities
\[
 sT_0=T_0s=e,\qquad sT_0s=s,\qquad T_0sT_0=T_0
\]
show that $s^\#=T_0$.
Moreover $t^2s=(T_0+\nu)^2s=T_0^2s=T_0$.
This agrees with~\cite[Theorem~3]{Drazin58}, applied to the
group-invertible element $T_0$ with inverse $s$; no finite-Drazin
invertibility of the full element $t$ is assumed here.
Since a group inverse is also the finite Drazin inverse,
$s^{\Dz}=T_0$, and this equals $t=T_0+\nu$ exactly when $\nu=0$.
The same arguments apply to the zero corners under our stated convention.
\end{proof}

\begin{definition}[Full and restricted Banach configurations]\label{def:classes}
In a complex unital Banach algebra, a full configuration consists of the coupled data of Notation~\ref{not:coupled}, with $pb=br$ and quasinilpotent input remainders. Coincidence follows from Theorem~\ref{thm:rigidity}. A restricted configuration is a full configuration whose residue $\nu$ is quasinilpotent in $B=(1-e)\A(1-e)$. A configuration is classical when $e=1$ and nondegenerate when $0\ne e\ne1$.
\end{definition}
The full class does not impose generalized-Drazin invertibility of $t$; even nilpotent input remainders do not imply that conclusion.
\begin{theorem}[Spectrum, generalized inverse, and the restricted resolvent]\label{thm:spectrum}
For a full Banach configuration,
\[
 \sigma(t)=\sigma_{e\A e}(T_0)\cup\sigma_B(\nu),\qquad t\in\A^{\gd}\ \Longleftrightarrow\ \nu\in B^{\gd}.
\]
When these inverses exist,
\[
 t^{\gd}=s+\nu_B^{\gd},\qquad 1-tt^{\gd}=f-\nu\nu_B^{\gd}.
\]
In particular, $s=t^{\gd}$ if and only if \(\nu\) is quasinilpotent in \(B\).

On the restricted class, $0$ is isolated in or absent from $\sigma(t)$, its Riesz idempotent is $f$, and for small nonzero $\lambda$
\begin{equation}\label{eq:resolvent}
 (\lambda1-t)^{-1}=(\lambda e-T_0)^{-1}_{e\A e}+\sum_{j=0}^\infty\lambda^{-j-1}\nu_B^j.
\end{equation}
Here $\nu_B^0=f$. If $f\ne0$, zero is a simple pole when $\nu=0$, a pole of order $k$ when $\nu$ is nilpotent of index $k$, and an isolated essential singularity otherwise. If $f=0$, the resolvent is holomorphic at zero.
\end{theorem}
\begin{proof}
Lemma~\ref{lem:blocks} proves the spectrum formula. The core spectrum is bounded away from zero, so zero is not an accumulation point of $\sigma(t)$ exactly when it is not an accumulation point of $\sigma_B(\nu)$. Koliha's spectral characterization~\cite[Theorem 4.2]{Koliha96} gives the inverse-existence equivalence. The element $s+\nu_B^{\gd}$ commutes with $t$, is an outer inverse, and has quasinilpotent remainder $\nu-\nu^2\nu_B^{\gd}$ in $B$. The same remainder is quasinilpotent in $\A$, by Lemma~\ref{lem:blocks}. Uniqueness~\cite[Theorem 4.3]{Koliha96} proves the formula and its spectral idempotent. The candidate equals this inverse exactly when $\nu_B^{\gd}=0$, equivalently when $\nu$ is quasinilpotent.

On the restricted class, the root test gives norm convergence of the series in~\eqref{eq:resolvent} for every $\lambda\ne0$, locally uniformly away from zero. Multiplication by $\lambda f-\nu$ telescopes to $f$, so it is the corner resolvent. Since $T_0$ is invertible in $e\A e$, its corner resolvent is holomorphic at zero. The coefficient of $\lambda^{-1}$ is $f$, giving the Riesz idempotent by integration on a small positively oriented circle. The remaining claims follow from the uniqueness of the Laurent coefficients: the principal part terminates precisely when some power of $\nu$ vanishes, and its highest nonzero term has order equal to the nilpotency index.
\end{proof}

On the restricted class, the analytic residue at zero of
$(\lambda1-t)^{-1}$ is $f$, the coefficient of $\lambda^{-1}$ in
\eqref{eq:resolvent}; it is zero when $f=0$.
This differs from our use of \emph{residue} for the residual parallel
sum $\nu$ of the Hua element.
\begin{remark}[The additional full-class cases]\label{rem:full}
Outside the restricted class, the output may still be generalized Drazin invertible, but its inverse has a nonzero residual term and its zero-spectral idempotent need not equal $f$. Alternatively, zero can accumulate in the residual spectrum and the output has no generalized Drazin inverse. The cases in which the candidate $s$ is the output inverse, therefore, do not exhaust the full class: the output can have a different generalized Drazin inverse, or it can have none. The general spectral machinery in Theorem~\ref{thm:spectrum} is classical; the coupling supplies the common core and the specific residual parallel sum to which it applies.
\end{remark}
\begin{proposition}[Functional calculus on the restricted class]\label{prop:calculus}
Let a restricted configuration be given, with $f\ne0$, and let $h$ be holomorphic on a neighborhood of $\sigma(t)$. If $h$ has no zeros on $\sigma(t)$, then $h(t)$ is invertible. Suppose instead that $h(0)=0$, that its order of vanishing at zero is the finite positive integer $m$, and that $h$ has no zero on $\sigma_{e\A e}(T_0)$. Then $h(t)$ has zero-spectral idempotent $f$ and supported generalized inverse $h(T_0)^{-1}_{e\A e}$. If $h(z)=z^mg(z)$ near zero, its residual part is $\nu^mg_B(\nu)$, with $g_B(\nu)$ invertible in $B$. A nonzero nilpotent residue of index $k$ is transformed to a residue of index $\lceil k/m\rceil$; a quasinilpotent nonnilpotent residue remains nonnilpotent. If $h$ is identically zero near zero, its residual part is zero. If, in addition, $h$ has no zero on $\sigma_{e\A e}(T_0)$, then $h(t)$ is group invertible with group inverse $h(T_0)^{-1}_{e\A e}$ and zero-spectral idempotent $f$.
\end{proposition}
\begin{proof}
The first assertion follows from the spectral mapping theorem.
Since the resolvent of $t$ is block diagonal along $e$, the bounded
corner compressions pass through the holomorphic contour integral. Thus
\[
 h(t)=H_0+\mu,\qquad
 H_0=h_{e\A e}(T_0),\qquad \mu=h_B(\nu),
\]
where the two corner calculi use the identities $e$ and $f$, respectively.

Suppose that $h$ has a zero of finite positive order $m$ at zero and
has no zeros on $\sigma_{e\A e}(T_0)$. Spectral mapping gives
$H_0\in(e\A e)^\times$ and, since $h(z)=z^mg(z)$ near zero with
$g(0)\ne0$,
\[
 \mu=\nu^m g_B(\nu),\qquad
 g_B(\nu)\in B^\times,\qquad \sigma_B(\mu)=\{0\}.
\]
The inverse $y=(H_0)^{-1}_{e\A e}$, extended by zero on $B$,
is a commuting outer inverse of $h(t)$ with
$h(t)-h(t)^2y=\mu$ and $h(t)y=e$.
Hence $y=h(t)^{\gd}$ and $1-h(t)y=f$.
Since $g_B(\nu)$ commutes with $\nu$ and is invertible in $B$,
\[
 \mu^j=\nu^{mj}g_B(\nu)^j,\qquad
 \mu^j=0\ \Longleftrightarrow\ \nu^{mj}=0
 \quad(j\ge1).
\]
This gives the nilpotency index $\lceil k/m\rceil$ when $\nu$
has index $k$, and preserves nonnilpotence when $\nu$ is not nilpotent.

If $h$ vanishes identically near zero, then $\mu=0$.
If $h$ also has no zeros on $\sigma_{e\A e}(T_0)$, the inverse of
$H_0$ in $e\A e$, extended by zero on $B$, is the group inverse
of $h(t)$, with zero-spectral idempotent $f$.
\end{proof}

\begin{remark}\label{rem:zero-core-spectrum}
Retain the restricted-class hypotheses of
Proposition~\ref{prop:calculus}, including $f\ne0$ and a zero
of finite positive order for $h$ at $0$, but allow $h$ to
vanish on the core spectrum. Put
$H_0=h_{e\A e}(T_0)$. If
$h(\lambda_0)=0$ for some
$\lambda_0\in\sigma_{e\A e}(T_0)$, then
$0\in\sigma_{e\A e}(H_0)$, and there are two possibilities.
If zero is isolated in this core spectrum, let $q_h\ne0$ be
its Riesz idempotent in $e\A e$. Then $h(t)$ is generalized
Drazin invertible, with
\[
 h(t)^{\gd}=(H_0)^{\gd}_{e\A e},\qquad
 1-h(t)h(t)^{\gd}=f+q_h.
\]
Here $fq_h=q_hf=0$ and $q_h\ne0$, so this idempotent
differs from $f$. If zero is an accumulation point
of $\sigma_{e\A e}(H_0)$, then $h(t)$ has no generalized
Drazin inverse and no Riesz idempotent at zero. These statements
follow from the two-corner spectral decomposition and the
quasinilpotence of $h_B(\nu)$.

The isolated alternative occurs whenever the core spectrum is
finite and $h$ has a zero on it. For an accumulating example,
let $\iota(x)=x$ in $C([1,2])$. The classical data
\[
 a=c=2\iota,\qquad b=(4\iota)^{-1},\qquad
 a'=c'=(2\iota)^{-1}
\]
have Hua element $\iota$. Adjoin the configuration of
Example~\ref{ex:nilpotent}. In
$C([1,2])\oplus\C\oplus M_3(\C)$ the resulting output is
\[
 t=\iota\oplus\tfrac12\oplus\tfrac12E_{12},
 \qquad
 \sigma_{e\A e}(T_0)=[1,2]\cup\{\tfrac12\}.
\]
For $h(z)=z(z-1)$ the zero at $0$ is simple, while
$h([1,2])=[0,2]$ accumulates at zero. Thus the transformed
output is not generalized Drazin invertible.

The holomorphic domain may be disconnected. On disjoint
neighborhoods of $[1,2]$, $0$, and $\tfrac12$, instead choose
$\widetilde h$ to be zero near $[1,2]$ and equal to $z$ near
$0$ and $\tfrac12$. Then zero is isolated, since
\[
 \widetilde h(t)=0\oplus\tfrac12\oplus\tfrac12E_{12},
 \qquad \sigma(\widetilde h(t))=\{0,\tfrac12\}.
\]

If instead $h(0)\ne0$, the residual block $h_B(\nu)$ is
invertible, so a zero-spectral idempotent of $h(t)$, when it
exists, need not contain $f$. The identically-zero germ at
zero is the separate case treated in
Proposition~\ref{prop:calculus}. None of these functional-calculus
statements constructs new coupled inputs realizing $h(t)$.
\end{remark}

\section{Commuting remainders and classicality}\label{sec:commuting}

\begin{lemma}[Quasinilpotent products, corners, and units]\label{lem:qnil}
Let $\R$ be a unital ring.
\begin{enumerate}
\item[(a)] A unit in a nonzero ring is not Harte-quasinilpotent.
If $q$ is Harte-quasinilpotent and $v,x\in\comm(q)$, then
$1+qvx$ is a unit.
\item[(b)] If $q$ is Harte-quasinilpotent and commutes with an
idempotent $f$, then $fqf$ is Harte-quasinilpotent in $f\R f$,
whose unit is $f$.
\item[(c)] If $q_1,q_2$ are commuting Harte-quasinilpotents, then
\[
1+(q_1+q_2)y\in\R^\times
\qquad\bigl(y\in\comm(q_1)\cap\comm(q_2)\bigr).
\]
Consequently $q_1+q_2$ is Harte-quasinilpotent if every element
commuting with $q_1+q_2$ also commutes with $q_1$.
\item[(d)] If $u\in\R^\times$ and $q$ is a
Harte-quasinilpotent commuting with $u$, then $u+q$ is a unit.
\end{enumerate}
\end{lemma}
\begin{proof}
For (a), a quasinilpotent unit would make $1-q q^{-1}=0$
a unit; the product assertion follows from $vx\in\comm(q)$.
For (b), if $y\in f\R f$ commutes with $fqf$, then, viewed in
$\R$, it commutes with $q$. Thus $1+qy$ is invertible and
commutes with $f$; its $f$-corner $f+(fqf)y$ is invertible
in $f\R f$.
For (c), put $u=1+q_1y$. This is a unit, and $q_2$ commutes
with $u$ and $u^{-1}$. Since $u^{-1}y\in\comm(q_2)$,
\[
1+(q_1+q_2)y=u\bigl(1+q_2u^{-1}y\bigr)
\]
is a product of units. If $y$ commutes with $q_1+q_2$ and
with $q_1$, it also commutes with $q_2$, proving the conditional
consequence. Finally, (d) follows from
$u+q=u(1+qu^{-1})$.
\end{proof}
\begin{lemma}[A commuting quasinilpotent sum cannot be a unit]\label{lem:unit}
If $q_1,q_2$ are commuting Harte-quasinilpotents in a nonzero unital ring, then $q_1+q_2$ is not invertible.
\end{lemma}
\begin{proof}
Suppose $u=q_1+q_2$ were invertible. Its inverse commutes with both summands. Quasinilpotence of $q_1$ makes $1-q_1u^{-1}=q_2u^{-1}$ invertible. Thus $q_2$ is a unit, contrary to Lemma~\ref{lem:qnil}.
\end{proof}
The argument is a short consequence of the established quasinilpotence definition and invertibility identities, not an independent generalized-inverse theory. For more details, see~\cite[formula (2.1)]{KolihaPatricio02} and~\cite[Lemma 2.3]{CuiChen13}.
\begin{corollary}\label{cor:idempotent}
A nonzero idempotent is not a sum of two commuting Harte-quasinilpotents.
\end{corollary}
\begin{proof}
If $e=q_1+q_2=e^2\ne0$, then both $q_i$ commute with $e$. Their compressions in the nonzero ring $e\R e$ are quasinilpotent by Lemma~\ref{lem:qnil}, commute, and sum to the unit $e$, contradicting Lemma~\ref{lem:unit}.
\end{proof}
\begin{proposition}[Commuting input remainders]\label{prop:commuting}
Assume coincidence and Harte-quasinilpotent input remainders in a unital ring. If $n_an_c=n_cn_a$, then $e=1$ and the configuration is classical. In particular, if either input remainder is zero, the configuration is classical. Thus a nonclassical intertwined finite-Drazin configuration has both input indices at least two.
\end{proposition}
\begin{proof}
In $f\R f$ the two remainders are quasinilpotent, commute, and have invertible sum with inverse $fbf$. Lemma~\ref{lem:unit} forces this corner to be the zero ring, so $f=0$. Both inputs are then units by Lemma~\ref{lem:input}. A zero remainder commutes with the other one; group-inverse input data have zero remainder.
\end{proof}

For commuting input remainders of a common-core configuration in a normed algebra, vanishing of both root-norm radii also forces classicality. Indeed, for any $\epsilon>0$ there are constants $K_a,K_c$ with $\|n_a^k\|\le K_a\epsilon^k$ and $\|n_c^k\|\le K_c\epsilon^k$ for every $k\ge0$. The binomial formula gives
\[
 \|(n_a+n_c)^k\|\le K_aK_c(2\epsilon)^k.
\]
Put $u=n_a+n_c$ and let $w=fbf$ be its inverse in the complementary corner, whose unit is $f=1-e$. If $f\ne0$, then for $k\ge1$,
\[
 0<\|f\|=\|w^ku^k\|\le\|w\|^k\|u^k\|.
\]
Taking $k$th roots gives
\[
 \liminf_{k\to\infty}\|u^k\|^{1/k}\ge\|w\|^{-1}>0,
\]
contrary to the binomial estimate. Thus $f=0$ and the configuration is classical.
\begin{proposition}[Normal input]\label{prop:normal}
In a full Banach configuration in a unital $C^*$-algebra, if either input is normal, then $e=1$.
\end{proposition}
\begin{proof}
Suppose $a$ is normal. Its generalized Drazin inverse and core idempotent belong to the commutative $C^*$-algebra generated by $a$ and $1$, by the holomorphic functional calculus. Thus $a(1-e)$ is normal and quasinilpotent, hence zero. Proposition~\ref{prop:commuting} applies. The argument for $c$ is identical.
\end{proof}
The hypotheses include quasinilpotent input remainders. A normal element with an arbitrary selected outer inverse is a different setting.
\begin{example}[Commuting nonclosure in a ring]\label{ex:nonclosure}
For this example we leave the Banach setting and work in the unital ring
\[
 R=\begin{pmatrix}\mathbb Z_{(\ell)}&\mathbb Q\\0&\mathbb Q\end{pmatrix},\qquad
 q_1=\ell I_2+E_{12},\qquad q_2=-E_{12},
\]
where $\ell$ is a prime. The two elements commute, and $q_2^2=0$. Any $X\in R$ commuting with $q_1$ must commute with $E_{12}$, hence has the form $X=dI_2+uE_{12}$ with $d\in\mathbb Z_{(\ell)}$ and $u\in\mathbb Q$. Then
\[
 1+q_1X=(1+\ell d)I_2+(\ell u+d)E_{12}
\]
is a unit of $R$, because $1+\ell d$ is a unit of the local ring. Thus $q_1$ and $q_2$ are Harte-quasinilpotent. Their sum $q_1+q_2=\ell I_2$ is central, but for $Y=-\ell^{-1}E_{22}\in R$,
\[
 1+(q_1+q_2)Y=E_{11}\notin R^\times.
\]
So the sum is not Harte-quasinilpotent. The detecting $Y$ does not belong to the joint commutant of the summands; the example, therefore, respects Lemma~\ref{lem:qnil}(c) and the unit-avoidance statements. Nonclosure itself was already recorded by Koliha and Patr\'icio~\cite[page 138]{KolihaPatricio02}. This explicit witness shows that one of the two commuting summands can be square-zero. It makes no assertion about sums of spectral-radius-zero elements in a complex Banach algebra.
\end{example}

\section{Examples separating the output regimes}\label{sec:examples}
The following examples establish the distinctions used here. They are not a fixed-dimensional classification of all residual parallel sums.

\emph{The square-zero threshold.}
For common-core coupled data, put $e=p=r$, $f=1-e$ and
$\beta=(n_a+n_c)^{-1}_{f\R f}$. If $n_a^2=0$, then
$\nu=n_a-n_a\beta n_a$ gives $\nu^2=0$; the symmetric formula
treats $n_c^2=0$. By Theorem~\ref{thm:output}, either condition
therefore gives $s=t^{\Dz}$ and $\ind_{\R}(t)\le2$.
If both squares vanish, then
$n_c(n_a+n_c)=(n_a+n_c)n_a$, so $n_a\beta=\beta n_c$ and
$\nu=\beta n_c^2=0$; hence $s=t^\#$.
Thus a nonnilpotent residual parallel sum from nilpotent input
remainders requires both input nilpotency indices to be at least
three. Examples~\ref{ex:family} and~\ref{ex:infinite} attain this
threshold.
\begin{example}[Vanishing residue]\label{ex:zero}
In $M_3(\C)$ put $a=E_{11}+E_{23}$, $c=E_{11}+E_{32}$, and $a'=c'=E_{11}$. Then
\[
 b=\tfrac12E_{11}+E_{23}+E_{32},\quad e=E_{11},\quad n_a=E_{23},\quad n_c=E_{32},\quad t=\tfrac12E_{11},\quad s=2E_{11}.
\]
Both remainders are nonzero and square-zero, their sum is invertible in the complementary corner, and their residual parallel sum is zero. These equations verify all the hypotheses directly. The configuration is nondegenerate and nonclassical, although $s=t^\#$.
\end{example}
\begin{example}[A nonzero nilpotent residue]\label{ex:nilpotent}
On a three-dimensional complementary corner put
\[
 N=E_{12}+E_{23},\quad C=E_{12}+E_{31},\quad
 \beta=\begin{pmatrix}0&0&1\\1/2&0&0\\0&1&0\end{pmatrix}.
\]
Then $N^3=C^3=0$, $(N+C)\beta=\beta(N+C)=I$, and $N\beta C=E_{12}/2$. Adjoining a scalar core gives, in $\C\oplus M_3(\C)$,
\[
 a=1\oplus N,\quad c=1\oplus C,\quad a'=c'=1\oplus0,\quad b=\tfrac12\oplus\beta,\quad t=\tfrac12\oplus E_{12}/2.
\]
Thus $s=2\oplus0=t^{\Dz}$, the output has index two, and $s$ is not its group inverse. The same example can be realized by block-diagonal matrices in $M_4(\C)$.
\end{example}
\begin{example}[A one-parameter residue family]\label{ex:family}
In $M_4(\C)$ let $e=E_{11}$, put $f=1-e$, and on the complementary corner carried by indices $2,3,4$ let
\[
 n_a=N:=E_{23}+E_{34},\qquad
 n_c=C_x:=-xE_{22}+E_{23}-x^2E_{32}+xE_{33}+E_{42}
 \qquad(x\in\C).
\]
Here $C_x=SC_0S^{-1}$ for the unit $S=1+xE_{32}$, so both input remainders are nilpotent of index three for every $x$. In the ordered complementary basis $(2,3,4)$, put
\[
 N+C_x=\begin{pmatrix}-x&2&0\\-x^2&x&1\\1&0&0\end{pmatrix},
 \qquad
 \beta_x=\begin{pmatrix}
 0&0&1\\[1pt]
 \tfrac12&0&\tfrac x2\\[1pt]
 -\tfrac x2&1&\tfrac{x^2}2
 \end{pmatrix}.
\]

Direct multiplication gives $(N+C_x)\beta_x=\beta_x(N+C_x)=I_3$, where $I_3$ represents the corner unit $f$, and
\[
 N\beta_xC_x=
 \begin{pmatrix}0&\tfrac12&0\\0&\tfrac x2&0\\0&0&0\end{pmatrix}
 =R_x:=\tfrac12E_{23}+\tfrac x2E_{33}.
\]
Complete the coupled data in $M_4(\C)$ by
\[
 a=e+N,\quad c=e+C_x,\quad a'=c'=e,\quad
 b=\tfrac12e+\beta_x,\quad b^{-1}=2e+N+C_x,
\]
where $\beta_x$ is embedded in $fM_4(\C)f$. The outer-inverse equations and $p=r=e$ follow directly. The element $b$ commutes with $e$, so intertwining holds. These are full configurations of Definition~\ref{def:classes}, with $T_0=\tfrac12e$, $s=2e$, and
\[
 \nu=R_x,\qquad R_x^2=\tfrac x2R_x,\qquad
 \sigma_{fM_4(\C)f}(R_x)=\{0,\tfrac x2\}.
\]
At $x=0$ the residue is nonzero and square-zero, so the configuration is restricted and the output has finite Drazin index two. For $x\ne0$ the residue carries the nonzero eigenvalue $x/2$, and the configuration leaves the restricted class although both input remainders are nilpotent. At $x=2$ the residue $R_2$ is idempotent, and Theorem~\ref{thm:spectrum} gives
\[
 t^{\gd}=2e+R_2\ne s=2e,\qquad
 1-tt^{\gd}=(1-e)-R_2.
\]
Along $x_n=1/n$ the raw and augmented data converge with $e_n=e$ and $s_n=s$ constant, while the eigenvalue $1/(2n)$ of $t_{x_n}$ tends to zero. Thus a uniform output gap fails on the full class even when $e_n$ and $s_n$ are constant, as used in Theorem~\ref{thm:full-continuity}.
\end{example}
\begin{example}[Two infinite-dimensional boundary choices]\label{ex:infinite}
The two constructions use the same matrix pattern: a coefficient
operator $Q$ controls the nonzero residual spectral values, up to a
factor of $1/2$, while the pattern adds a zero spectral value.
The choices below separate quasinilpotent nonnilpotence from
accumulation of nonzero spectrum at zero. Let $\mathcal K=B(C[0,1])$ and let $Q$ be the Volterra operator, or let $\mathcal K=B(\ell^2)$ and let $Q=\diag(2,1,2/3,\ldots)$. With $1$ denoting the identity of $\mathcal K$, put
\[
 N=\begin{pmatrix}0&1&0\\0&0&1\\0&0&0\end{pmatrix},\quad
 C_Q=\begin{pmatrix}-Q&1&0\\-Q^2&Q&0\\1&0&0\end{pmatrix},\quad
 \beta_Q=\begin{pmatrix}0&0&1\\1/2&0&Q/2\\-Q/2&1&Q^2/2\end{pmatrix}.
\]
Direct multiplication gives $N^3=C_Q^3=0$, $\beta_Q=(N+C_Q)^{-1}$, and
\[
 R_Q=N\beta_QC_Q=\tfrac12\begin{pmatrix}0&1&0\\0&Q&0\\0&0&0\end{pmatrix},\quad
 R_Q^j=2^{-j}\begin{pmatrix}0&Q^{j-1}&0\\0&Q^j&0\\0&0&0\end{pmatrix}\qquad(j\ge1).
\]
The scalar-core construction of Example~\ref{ex:nilpotent}, therefore, gives full configurations in $\C\oplus M_3(\mathcal K)$. For $\lambda\ne0$, direct multiplication gives
\[
 \lambda I_3-R_Q=
 \begin{pmatrix}
 \lambda1&0&0\\0&\lambda1-Q/2&0\\0&0&\lambda1
 \end{pmatrix}
 \begin{pmatrix}
 1&-(2\lambda)^{-1}1&0\\0&1&0\\0&0&1
 \end{pmatrix}.
\]

The second factor is a unit, with the sign of its off-diagonal entry reversed in its inverse. The diagonal factor is invertible exactly when $\lambda1-Q/2$ is invertible in $\mathcal K$. Also $R_QE_{11}=0$ with $E_{11}\ne0$, so $R_Q$ is not invertible. Consequently
\[
 \sigma_{M_3(\mathcal K)}(R_Q)=\{0\}\cup\tfrac12\sigma_{\mathcal K}(Q).
\]

For Volterra, $\|Q^j\|=1/j!$ and no power vanishes, as follows from the iterated integration formula and its value on the constant function $1$. The displayed powers show that $R_Q$ is quasinilpotent and nonnilpotent. The configuration is restricted, but its output has no finite Drazin inverse. For the diagonal choice, $\sigma(R_Q)=\{0\}\cup\{1/n:n\ge1\}$, so the output is not generalized Drazin invertible. Both conclusions hold although the input remainders are nilpotent of index three. Only these concrete boundary instances are needed here; the general realization problem is separate.
\end{example}
\begin{example}[A strict ring generalized inverse without finite Drazin index]\label{ex:product}
Let $\R=\mathbb Q\times\prod_{d\ge3}M_d(\mathbb Q)$. In a direct product, commutants, units and idempotents are all computed componentwise; consequently a tuple whose components are nilpotent of unbounded index is Harte-quasinilpotent without being nilpotent. In the factor $M_d(\mathbb Q)$, with standard basis $f_1,\ldots,f_d$, define nilpotent matrices
\[
 n_a^{(d)}:\ f_j\mapsto f_{j-1}\ (2\le j\le d);\qquad
 n_c^{(d)}:\ f_j\mapsto f_{j-1}\ (2\le j\le d-1),\quad f_1\mapsto f_d.
\]
All basis vectors not listed are mapped to zero. Their sum $S^{(d)}$ maps $f_j\mapsto2f_{j-1}$ $(2\le j\le d-1)$, $f_d\mapsto f_{d-1}$ and $f_1\mapsto f_d$, an invertible weighted cycle whose inverse $b_d$ acts by $f_i\mapsto\tfrac12f_{i+1}$ $(1\le i\le d-2)$, $f_{d-1}\mapsto f_d$, $f_d\mapsto f_1$. A basis chase gives
\[
 \nu_d:=n_a^{(d)}b_dn_c^{(d)}:\ f_j\mapsto\tfrac12f_{j-1}\ (2\le j\le d-1),\qquad f_1,f_d\mapsto0,
\]
nilpotent of index exactly $d-1$; for $d=3$ this reproduces the residue in Example~\ref{ex:nilpotent}. Now fix $\alpha,\gamma\in\mathbb Q^\times$ with $\alpha+\gamma\ne0$ and set
\[
 a=\bigl(\alpha,(n_a^{(d)})_{d\ge3}\bigr),\qquad
 c=\bigl(\gamma,(n_c^{(d)})_{d\ge3}\bigr),\qquad b^{-1}=a+c\in\R^\times.
\]
Then $a'=(\alpha^{-1},0,0,\ldots)$ and $c'=(\gamma^{-1},0,0,\ldots)$ are commuting outer inverses with quasinilpotent input remainders and common idempotent $e=(1,0,0,\ldots)$. Neither input is finite-Drazin invertible. Indeed, if $n'$ is a commuting outer inverse of a nilpotent $n$ and $p=nn'$, then $np$ is both nilpotent and invertible in $p\R p$. Thus $p=0$ and $n'=n'p=0$. Applied componentwise, this forces zero inverse components on the matrix factors, leaving input remainders with unbounded nilpotency indices.

Theorem~\ref{thm:common} applies, and the residue of $t=a-aba$ is $\nu=(0,\nu_3,\nu_4,\ldots)$, which is Harte-quasinilpotent but not nilpotent. The element $t$ is not finite-Drazin invertible: the axiom $t^k=t^{k+1}t^{\Dz}$ would give, in block $d$, $\nu_d^k=\nu_d^{k+j}(t_d^{\Dz})^j=0$ for sufficiently large $j$, forcing $k\ge d-1$ for every $d$. On the other hand, $s=a'+c'$ is central in $\R$, since only its scalar component $\alpha^{-1}+\gamma^{-1}$ is nonzero. Hence $s\in\comm^2(t)$, and
\[
 sts=s,\qquad ts=st=e,\qquad t-t^2s=\nu\in\qnil(\R).
\]
Thus $s$ is the strict ring generalized Drazin inverse of $t$ in the sense of~\cite{KolihaPatricio02}, while $t$ has no finite Drazin inverse.

View the same residue in the complex Banach algebra of bounded tuples
\[
 \C\times\prod_{d\ge3}^{\mathrm{bdd}}M_d(\C),
\]
with the supremum of the Euclidean operator norms. In this algebra $\rho(\nu)=\tfrac12$. Indeed, for every $k\ge1$, the matrices $\nu_d^k$ are $2^{-k}$-weighted shifts and have norm $2^{-k}$ whenever $d\ge k+2$, so $\|\nu^k\|=2^{-k}$. Consequently $\nu$ is not Harte-quasinilpotent in this bounded Banach algebra. The ring--Banach dictionary compares the two definitions within one algebra; changing the ambient algebra changes the relevant unit tests. A Banach configuration with quasinilpotent nonnilpotent residue is given by the Volterra choice in Example~\ref{ex:infinite}.
\end{example}

\section{Canonicity and the two continuity criteria}\label{sec:continuity}
All algebras in this section are complex unital Banach algebras. Limits stated to be configurations include all hypotheses of Definition~\ref{def:classes}; they are not arbitrary limits of raw inputs.
\begin{lemma}[Compact-parameter inversion]\label{lem:compact}
Let $K$ be compact, let $x:K\to\A^\times$ be continuous, and suppose continuous $x_n:K\to\A$ converge uniformly to $x$. Then, for all sufficiently large $n$, $x_n$ is invertible on $K$ and $x_n^{-1}\to x^{-1}$ uniformly.
\end{lemma}
\begin{proof}
Let $M=\sup_K\|x^{-1}\|<\infty$. For large $n$, $\sup_K\|x^{-1}(x_n-x)\|<1/2$. The Neumann series gives the inverses and a uniform bound $2M$. The inverse-difference identity then gives uniform convergence.
\end{proof}
\begin{proposition}[Canonical inverse data and recovery]\label{prop:canonical}
In a full configuration, $a'=a^{\gd}$ and $c'=c^{\gd}$. The data $e,s$ are, therefore, determined by $(a,b)$. Moreover
\begin{equation}\label{eq:recover}
 a'=bcs=scb,\qquad c'=bas=sab,\qquad s=e(1-e+te)^{-1}.
\end{equation}
On the restricted class $s=t^{\gd}$ and $1-e$ is the zero-spectral idempotent of $t$.
\end{proposition}
\begin{proof}
The chosen witnesses satisfy the defining Banach generalized-Drazin equations, so uniqueness~\cite[Theorem 4.3]{Koliha96} proves the first assertion. This canonicity has the historical antecedent~\cite[Definition 7.5.2 and Theorem 7.5.3, pp.~257--259]{Harte88}. Harte's support is our core idempotent, not its complement. For recovery, take $bcs=bc(a'+c')=(1-ba)a'+be=a'$; likewise $scb=a'(1-ab)+eb=a'$. Interchanging $a,c$ gives the other two identities. The formula for $s$ is Theorem~\ref{thm:common}, and the restricted conclusions follow from Theorem~\ref{thm:spectrum}.
\end{proof}
For any $x\in\A$, define
\[
 \gap(x)=d(0,\sigma(x)\setminus\{0\}),\qquad d(0,\varnothing)=+\infty.
\]
If $x\in\A^{\gd}$, the holomorphic definition of its inverse and spectral mapping give
\begin{equation}\label{eq:gap}
 \gap(x)=1/\rho(x^{\gd}),\qquad 1/0=+\infty.
\end{equation}
For $\rho(x)>0$ this is~\cite[Lemma 1.3]{KolihaRakocevic98}. If the nonzero spectrum is empty, $x^{\gd}=0$; otherwise its nonzero spectral values are precisely the reciprocals of those of $x$. The norm obeys $\|x^{\gd}\|\ge1/\gap(x)$, with no general reverse bound.
\begin{theorem}[Continuity on the full class]\label{thm:full-continuity}
Let full configurations have $(a_n,b_n)\to(a,b)$ in norm, with a full limiting configuration. Then the following are equivalent:
\begin{enumerate}
\item[(i)] $e_n\to e$;
\item[(ii)] $s_n\to s$;
\item[(iii)] $\sup_n\|s_n\|<\infty$;
\item[(iv)] $\inf_n\gap(a_n)>0$;
\item[(v)] $\inf_n\gap(c_n)>0$.
\end{enumerate}
A uniform output gap is not equivalent to these conditions.
\end{theorem}
\begin{proof}
Continuity of inversion gives $c_n=b_n^{-1}-a_n\to c$ and multiplication gives $t_n\to t$. If (i) holds, the last formula of~\eqref{eq:recover} and continuity of inversion imply (ii). Clearly (ii) implies (iii). Under (iii), $a_n'=b_nc_ns_n$ is bounded, because both $b_n,c_n$ are bounded. Equation~\eqref{eq:gap}, applied to $x=a_n$ with $a_n^{\gd}=a_n'$, gives (iv).

For (iv)$\Rightarrow$(i), choose $0<\epsilon<\min\{\inf_n\gap(a_n),\gap(a)\}$. The limit is generalized Drazin invertible, so this is possible. The circle $|\lambda|=\epsilon$ lies in the resolvent sets of every $a_n$ and of $a$ and encloses none of their nonzero spectral values. Lemma~\ref{lem:compact} gives uniform resolvent convergence on the circle. Integration gives convergence of the zero-spectral idempotents $1-e_n\to1-e$. Using $c_n'=b_na_ns_n$ and the same argument proves equivalence with (v). In Example~\ref{ex:family}, the parameter $x_n=1/n$ leaves $e,s$ constant while an output eigenvalue $1/(2n)$ tends to zero. Thus the output gap is not controlled by these equivalent conditions.
\end{proof}
\begin{theorem}[The output criterion on the restricted class]\label{thm:restricted-continuity}
Under the hypotheses of Theorem~\ref{thm:full-continuity},
assume in addition that every configuration in the sequence
and the limiting configuration belong to the restricted class.
Then conditions \textup{(i)--(v)} of that theorem are also
equivalent to each of the following:
\begin{enumerate}
\item[(vi)] $\inf_n\gap(t_n)>0$;
\item[(vii)] for some $\epsilon>0$, every $t_n$ has no nonzero spectrum in $|\lambda|\le\epsilon$ and $\sup_n\sup_{|\lambda|=\epsilon}\|(\lambda1-t_n)^{-1}\|<\infty$.
\end{enumerate}
On this class $\gap(t)=1/\rho(s)$. Failure of $e_n\to e$ implies $\inf_n\gap(t_n)=0$ and $\sup_n\rho(s_n)=\infty$.
\end{theorem}
\begin{proof}
Here $s_n=t_n^{\gd}$ and $s=t^{\gd}$. Thus boundedness of $s_n$ implies (vi) by~\eqref{eq:gap}.

Conversely, (vi), together with generalized-Drazin invertibility of the limit, permits a common small circle enclosing only zero spectral parts. Lemma~\ref{lem:compact} and integration give $1-e_n\to1-e$. The full-class theorem completes the equivalence. The same circle has a uniformly bounded family of resolvents: eventual uniform convergence gives a bound on the tail, and the finite initial set has a bound by compactness. Hence (vi) implies (vii); (vii) itself implies (vi). The last assertions follow from~\eqref{eq:gap} and the equivalence just proved.
\end{proof}
\begin{remark}[Exact scope of the continuity sources]\label{rem:continuity-sources}
For a sequence of generalized-Drazin invertible elements with a generalized-Drazin invertible limit, Koliha and Rako\v{c}evi\'c~\cite[Theorem 2.4]{KolihaRakocevic98} establish the equivalence between convergence of the inverses, convergence of the spectral idempotents, boundedness of the inverses, and a positive uniform spectral gap. The restricted theorem is its application through canonicity. The full-class result instead uses the input inverses and the recovery identities, since the output may not be generalized Drazin invertible. We use no additional reduced-minimum-modulus criterion. For bounded operators of finite Drazin index the relevant inverse and core-product conditions also occur in~\cite[Theorem 2.2, conditions (2.2)--(2.4)]{Rakocevic99}. The qualitative closed-operator graph-gap theorem~\cite[Theorem 4.1]{CGKR02} is a separate statement with its own limiting-class and domain hypotheses, not a numerical estimate invoked here.
\end{remark}
\begin{example}[Failure of continuity and the inverse norm]\label{ex:continuity}
Use the limiting configuration of Example~\ref{ex:zero} and fix its $b$. Put
\[
 a_n=E_{11}+E_{23}+n^{-1}(E_{22}+E_{33}),\qquad c_n=b^{-1}-a_n.
\]
Both inputs are invertible, so $e_n=I$ and the configurations are classical, but $e=E_{11}$. Direct block inversion gives
\[
 s_n=2E_{11}-n^2(E_{23}+E_{32}),\qquad t_n=\tfrac12E_{11}-n^{-2}(E_{23}+E_{32}).
\]
For $n\ge2$, $\|s_n\|=\rho(s_n)=n^2$ and $\gap(t_n)=n^{-2}$ in the Euclidean operator norm. Thus the reciprocal-gap lower bound is attained.

For a contrasting family keep $a_n$ and put
\[
 \widetilde c_n=E_{11}+E_{32}-n^{-2}(E_{22}+E_{33}),\qquad b_n=(a_n+\widetilde c_n)^{-1}.
\]
The sum is invertible for every positive integer $n$: its complementary determinant is $n^{-2}(1-1/n)^2-1\ne0$. The pairs converge to the same limit and remain classical away from the limit. The complementary block of $s_n$ is
\[
 \begin{pmatrix}n-n^2&-n^2\\-n^4&n-n^2\end{pmatrix}.
\]
Its eigenvalues are $n-n^2\pm n^3$. Hence $\rho(s_n)=n^3+n^2-n$ for $n\ge2$, whereas $\|s_n\|\ge n^4$. Therefore, $\|s_n\|\gap(t_n)\to\infty$. Spectral gap and inverse norm are not interchangeable measures.
\end{example}

\section{Exact defects and local augmented-data estimates}\label{sec:defects}
We retain the coupled outer data $(a,b,c,a',c')$ of
Notation~\ref{not:coupled}, but do not assume either $pb=br$ or $p=r$.
We begin in an arbitrary unital associative ring, without imposing a
smallness condition on the input remainders, and derive exact identities
relating the difference $r-p$ to the intertwining defect $pb-br$.
Additional nilpotence hypotheses give terminating formulas; in the
Banach setting, quasinilpotence gives norm-convergent formulas, while
invertibility of the coefficient operators suffices for local estimates.
These estimates concern the augmented tuple, including the selected
inverses $a',c'$, rather than the raw pair $(a,b)$ alone.
At an exactly intertwined base with invertible coefficient operators,
we also determine the feasible first-order defect space and its sharp
norm constant.
\begin{definition}[Intertwining defect]\label{def:defect}
Set
\[
 \delta=pb-br,\qquad D=b^{-1}\delta b^{-1}=(a+c)p-r(a+c).
\]
The primed witnesses are part of the augmented tuple. Approximate intertwining means a small $\delta$ only when a norm is specified.
\end{definition}
\begin{lemma}[Defect reduction]\label{lem:defect}
With $z=n_cp$ and $v=rn_a$, the following identities hold:
\begin{align*}
 D&=c(p-r)+(p-r)a,\\
 rp&=p+n_ca'-Da',\qquad rp=r+c'n_a+c'D,\\
 z&=-n_cza'+n_cDa',\qquad v=-c'vn_a-c'Dn_a,\\
 r-p&=n_ca'-c'n_a-Da'-c'D.
\end{align*}
No nilpotence or quasinilpotence hypothesis is required.
\end{lemma}
\begin{proof}
Expand the definition of $D$. Multiplying $r(a+c)=(a+c)p-D$ on the right by $a'$ and on the left by $c'$ gives the two formulas for $rp$. Multiplication by $n_c$ on the left or by $n_a$ on the right gives the two recurrence formulas, using the corner annihilations in Lemma~\ref{lem:input}. Subtraction yields the final relation.
\end{proof}
\begin{proposition}[Finite tails]\label{prop:tails}
For every integer $M\ge1$,
\begin{align*}
 z&=(-1)^Mn_c^Mz(a')^M+\sum_{k=1}^M(-1)^{k-1}n_c^kD(a')^k,\\
 v&=(-1)^M(c')^Mvn_a^M+\sum_{k=1}^M(-1)^k(c')^kDn_a^k.
\end{align*}
Consequently
\begin{align*}
 r-p={}&-Da'-c'D+\sum_{k=1}^M(-1)^{k-1}\bigl(n_c^kD(a')^{k+1}+(c')^{k+1}Dn_a^k\bigr)\\
 &+(-1)^M\bigl(n_c^Mz(a')^{M+1}-(c')^{M+1}vn_a^M\bigr).
\end{align*}
These are finite identities in every ring.
\end{proposition}
\begin{proof}
Induction on $M$ substitutes the recurrence for the final occurrence of $z$ or $v$ in Lemma~\ref{lem:defect}. Multiply the resulting formulas by $a'$ on the right and by $c'$ on the left, and use $n_ca'=za'$ and $c'n_a=c'v$ in the final reduction identity.
\end{proof}
For $x\in\R$ denote left and right multiplication by $L_x$ and $R_x$. These are additive maps in a ring and bounded complex-linear maps in a Banach algebra.
\begin{theorem}[Closed defect representation]\label{thm:defect}
If $n_c^M=n_a^M=0$, then
\[
 r-p=-\sum_{k=0}^{M-1}(-1)^k\bigl(n_c^kD(a')^{k+1}+(c')^{k+1}Dn_a^k\bigr).
\]
In a complex Banach algebra with quasinilpotent input remainders, the same identity holds with the two norm-convergent infinite series. Equivalently, putting
\[
 U=I+L_{n_c}R_{a'},\qquad V=I+L_{c'}R_{n_a},
\]
these additive or bounded-linear operators are invertible in the respective setting and
\begin{equation}\label{eq:defect}
 r-p=\mathcal F(D),\qquad \mathcal F(D)=-U^{-1}(D)a'-c'V^{-1}(D).
\end{equation}
\end{theorem}
\begin{proof}
In the nilpotent case the tails and terminal $k=M$ summands in Proposition~\ref{prop:tails} vanish. The finite geometric inverses of $U,V$ give the displayed operator form. In the Banach case,
\[
 \|(L_{n_c}R_{a'})^k\|\le\|n_c^k\|\|(a')^k\|,\qquad
 \|(L_{c'}R_{n_a})^k\|\le\|(c')^k\|\|n_a^k\|.
\]
Both root limits are zero. The Neumann series converge in operator norm, and the tails in the finite formula tend to zero. This proves the infinite-series and operator forms. No estimate on an inverse is obtained by replacing its norm with a reciprocal eigenvalue gap.
\end{proof}
\begin{corollary}[A data-dependent estimate]\label{cor:estimate}
For the Banach hypotheses of Theorem~\ref{thm:defect},
\[
 \|r-p\|\le\bigl(\|U^{-1}\|\,\|a'\|+\|c'\|\,\|V^{-1}\|\bigr)\|b^{-1}\|^2\,\|\delta\|.
\]
There is also the explicit series bound obtained by replacing the two operator-inverse norms by $\sum_{k\ge0}\|n_c^k\|\|(a')^k\|$ and $\sum_{k\ge0}\|(c')^k\|\|n_a^k\|$, respectively.
\end{corollary}
\begin{proof}
Apply the norm to~\eqref{eq:defect} and use $D=b^{-1}\delta b^{-1}$. The preceding norm-convergent Neumann series give the alternative constants.
\end{proof}

\begin{proposition}[Local variation on augmented data]\label{prop:local}
Let $\A$ be a complex unital Banach algebra, and consider the
coupled outer-data tuples $(a,b,c,a',c')$ of
Notation~\ref{not:coupled}, without assuming intertwining.
Whenever the bounded operators $U,V$ are invertible,
\eqref{eq:defect} holds. At any fixed such tuple, the inverses
$U^{-1},V^{-1}$ and the coefficient map $\mathcal F$ vary
continuously with the augmented data, and
\[
 \|r-p\|\le
 \bigl(\|U^{-1}\|\,\|a'\|+\|c'\|\,\|V^{-1}\|\bigr)
 \|b^{-1}\|^2\,\|\delta\|
\]
has a uniform coefficient on a sufficiently small neighborhood
among these tuples. Quasinilpotence of the input remainders
is sufficient for the inverse conditions, but is not required
for this local assertion.

If an augmented tuple varies differentiably with a parameter
$\theta$, with $U_0,V_0$ invertible and $\delta_0=0$, then
\[
 \dot r_0-\dot p_0=
 \mathcal F_0\bigl(b_0^{-1}\dot\delta_0b_0^{-1}\bigr).
\]
Exact intertwining is required here only at the base tuple.
\end{proposition}
\begin{proof}
Write $z=n_cp$ and $v=rn_a$. The recurrence equations of
Lemma~\ref{lem:defect} give
\[
 Uz=(U-I)D,\qquad Vv=-(V-I)D.
\]
Thus
\[
 z=D-U^{-1}(D),\qquad v=V^{-1}(D)-D.
\]
Substitution into $r-p=za'-c'v-Da'-c'D$ gives
\[
 r-p=-U^{-1}(D)a'-c'V^{-1}(D)=\mathcal F(D).
\]
Taking norms and using $D=b^{-1}\delta b^{-1}$ proves the
estimate. The multiplication operators depend continuously
on the augmented coordinates. Invertibility is open in the
Banach algebra of bounded operators, and inversion there
is continuous, so the coefficients are locally bounded.
For a differentiable path, the coefficient map is differentiable
as well. Since $D_0=0$, differentiating
$r_\theta-p_\theta=\mathcal F_\theta(D_\theta)$ gives
$\dot r_0-\dot p_0=\mathcal F_0(\dot D_0)$.
The product rule and $\delta_0=0$ yield
$\dot D_0=b_0^{-1}\dot\delta_0b_0^{-1}$.
\end{proof}
\begin{remark}[No raw-data stability conclusion]\label{rem:raw}
The neighborhood above includes the selected inverses as coordinates. Near a singular limiting input, the unique inverses at nearby raw pairs can diverge, as Example~\ref{ex:continuity} shows. The proposition supplies neither a raw-input Lipschitz theorem, a transversality statement, nor a nearby exactly intertwined tuple realizing a prescribed repair. Those would require different hypotheses and additional arguments.
\end{remark}

The coefficient in Corollary~\ref{cor:estimate} bounds the norm of the
frozen map $\mathcal G(Y)=\mathcal F(b^{-1}Yb^{-1})$ by separate
multiplication and inverse norms. The sharp first-order constant instead
uses the restriction of this map to feasible defect directions.

\begin{theorem}[Exact feasible-tangent constant]\label{thm:tangent-sharp}
Let $(a,b,c,a',c')$ be coupled outer data in a complex unital Banach
algebra $\A$, with $\delta=pb-br=0$ and with the operators $U,V$ of
Theorem~\ref{thm:defect} invertible. Then $p=r=e$.
Put $f=1-e$ and equip
\[
 W=e\A f\oplus f\A e
\]
with the norm inherited from $\A$. The bounded linear map
\[
 \mathcal L:W\longrightarrow W,\qquad
 \mathcal L(X)=-b(cX+Xa)b
\]
is invertible, and
\[
 \mathcal L^{-1}=\mathcal G|_W,\qquad
 \mathcal G(Y)=\mathcal F(b^{-1}Yb^{-1}).
\]
For every norm-differentiable path of coupled outer-data tuples through
this base tuple, with derivatives taken at zero,
\[
 X=\dot r_0-\dot p_0\in W,\qquad
 \dot\delta_0=\mathcal L(X).
\]
Conversely, every $X\in W$ is realized by an analytic path of such
tuples. Intertwining is required only at the base, not along the path.
The least nonnegative constant valid for all these paths in
\[
 \|\dot r_0-\dot p_0\|\le
 \kappa_{\mathrm{tan}}\|\dot\delta_0\|
\]
is therefore
\[
 \kappa_{\mathrm{tan}}=\|\mathcal L^{-1}\|_{\mathcal B(W)}.
\]
Here $\mathcal B(W)$ denotes the bounded linear operators on $W$;
if $W=\{0\}$, the constant is zero.
\end{theorem}
\begin{proof}
Proposition~\ref{prop:local} gives $r-p=\mathcal F(0)=0$.
The bounded projection $Z\mapsto eZf+fZe$ has range $W$, so $W$
is a Banach space. Differentiating $p_\theta^2=p_\theta$ and
$r_\theta^2=r_\theta$ at $p_0=r_0=e$ makes both derivatives
off-diagonal. Moreover, Lemma~\ref{lem:defect} gives
\[
 D_\theta=c_\theta(p_\theta-r_\theta)
             +(p_\theta-r_\theta)a_\theta,\qquad
 \delta_\theta=b_\theta D_\theta b_\theta.
\]
Since $D_0=0$, differentiation proves
$\dot\delta_0=-b(cX+Xa)b$.

Write $a=A_0+n_a$ and $c=C_0+n_c$ along $e$.
For $u\in e\A f$ and $v\in f\A e$, the two corner restrictions of
$X\mapsto cX+Xa$ are
\[
 u\longmapsto C_0u+un_a=C_0V(u),\qquad
 v\longmapsto n_cv+vA_0=U(v)A_0.
\]
The operators $U,V$ commute with all four corner projections, so their
inverses preserve the same corners. The core factors $A_0,C_0$ are
units in $e\A e$ by Lemma~\ref{lem:input}, and $b,b^{-1}$ commute
with $e$ by Theorem~\ref{thm:common}. These factorizations show that
$\mathcal L$ is invertible on $W$. The same corner calculation gives
$\mathcal F(cX+Xa)=-X$; hence $\mathcal G\mathcal L(X)=X$, proving
the inverse formula.

For feasibility, fix $X\in W$ and put
\[
 H=eXf-fXe,\qquad w_\theta=\exp(\theta H),\qquad
 a_\theta=w_\theta a w_\theta^{-1},\quad
 a'_\theta=w_\theta a'w_\theta^{-1},
\]
\[
 c_\theta=c,\qquad c'_\theta=c',\qquad
 b_\theta=(a_\theta+c)^{-1}.
\]
Invertibility is open, so these analytic tuples are defined for small
$\theta$ and satisfy every coupled outer-data equation.
The operators $U_\theta,V_\theta$ remain invertible after shrinking
the parameter neighborhood. Here
$p_\theta=w_\theta e w_\theta^{-1}$ and $r_\theta=e$;
since $[H,e]=-X$, their derivative difference is $X$.
Also $n_{a,\theta}=w_\theta n_a w_\theta^{-1}$ and
$n_{c,\theta}=n_c$, so any nilpotence or quasinilpotence of the input
remainders is preserved.

If $W\ne\{0\}$, bijectivity and feasibility give
\[
 \sup_{0\ne X\in W}\frac{\|X\|}{\|\mathcal L(X)\|}
 =\|\mathcal L^{-1}\|_{\mathcal B(W)}.
\]
Thus no smaller constant works for all paths. The supremum need not
be attained in infinite dimensions. If $W=\{0\}$, both derivatives
vanish and the least nonnegative constant is zero.
\end{proof}

For $W\ne\{0\}$, the same constant can be written as a reciprocal
separation:
\[
 \operatorname{sep}_{\mathrm{Hua}}
 =\inf_{0\ne X\in W}\frac{\|b(cX+Xa)b\|}{\|X\|},
 \qquad \kappa_{\mathrm{tan}}=\operatorname{sep}_{\mathrm{Hua}}^{-1}.
\]
This is the lower bound of the already defined weighted Sylvester map
on $W$. The standard matrix separation in the Frobenius norm is the
least singular value of the corresponding Sylvester operator; see
\cite[``Computing $s$ and $\mathrm{sep}$'']{LAPACK99}.
That Frobenius-norm formula is not automatically the constant for the
inherited algebra norm used here, nor can separation generally be
replaced by an eigenvalue gap.

The norm on $W$ is inherited; a general idempotent $e$ need not be
self-adjoint. The theorem concerns the entire local class of augmented
coupled outer data, not paths constrained to remain exactly intertwined.
It also applies at an interior point of additional strict bounds.
A smaller constrained class may have fewer feasible tangent directions,
and its sharp constant requires the corresponding restricted
optimization. In particular, neither $\|\mathcal G\|$ on all of $\A$
nor the product bound of Corollary~\ref{cor:estimate} is asserted to be
the sharp tangent constant. No globally uniform nonlinear coefficient
or nearby exact repair is supplied.

\begin{example}[Euclidean-norm tangent calibration]\label{ex:tangent-calibration}
Use the data of Example~\ref{ex:zero} in $M_3(\C)$, equipped with the
Euclidean operator norm:
\[
 a=E_{11}+E_{23},\qquad c=E_{11}+E_{32},\qquad
 b=\tfrac12E_{11}+E_{23}+E_{32},\qquad a'=c'=e=E_{11}.
\]
The input remainders are square-zero, so $U,V$ are invertible.
Relative to the orthogonal decomposition $\C^3=e\C^3\oplus f\C^3$,
write
\[
 X=\begin{pmatrix}0&u\\v&0_2\end{pmatrix},\qquad
 u\in M_{1,2}(\C),\quad v\in M_{2,1}(\C).
\]
The tangent map is
\[
 \mathcal L(X)=-\begin{pmatrix}0&uK\\Kv&0_2\end{pmatrix},
 \qquad K=\frac12\begin{pmatrix}1&1\\1&0\end{pmatrix}.
\]
Since $X^*X$ is block diagonal,
$\|X\|=\max\{\|u\|_2,\|v\|_2\}$, and consequently
$\|\mathcal L^{-1}\|=\|K^{-1}\|_2$.
The eigenvalues of the real symmetric matrix $K$ are
$(1+\sqrt5)/4$ and $(1-\sqrt5)/4$. Its least singular value is
$(\sqrt5-1)/4$, giving
\[
 \boxed{\kappa_{\mathrm{tan}}=1+\sqrt5.}
\]
Equality is attained by $v=0$ and a nonzero row eigenvector $u$ for
$(1-\sqrt5)/4$, with the corresponding direction realized by the
path in Theorem~\ref{thm:tangent-sharp}. This is a derivative ratio
at the specified base tuple, not a nonlinear bound with the same
constant on a whole neighborhood.
\end{example}

\section{Remaining structural questions}\label{sec:questions}
The following questions distinguish the cases settled above from the
remaining classification, quantitative and compatibility problems.
The common-core inverse criteria, finite-central intertwining theorem,
Noetherian power separation and pointwise tangent theorem provide
starting cases, not solutions of all the questions below.

\begin{problem}[Bare Drazin equalities]\label{prob:bare-drazin}
Let $\R$ be a unital associative ring, and let $a,c\in\R$ be
finite-Drazin invertible with $a+c\in\R^\times$.
Put $b=(a+c)^{-1}$ and $t=a-aba$.
Characterize the pairs for which $t$ is finite-Drazin invertible and
\[
 t^{\Dz}=a^{\Dz}+c^{\Dz},
\]
without assuming coincidence of the input core idempotents
$aa^{\Dz}$ and $cc^{\Dz}$.
\end{problem}

Under coincidence, Theorem~\ref{thm:output}(c) characterizes the equality
by nilpotence of the residue. Problem~\ref{prob:bare-drazin} asks for a
structural characterization without that assumption; it does not
impose the prescribed product identities of
Theorem~\ref{thm:common}.

\begin{problem}[Structural conditions for coincidence]\label{prob:structural-coincidence}
For coupled commuting outer data with Harte-quasinilpotent input
remainders and $pb=br$, identify natural ring or bimodule hypotheses
that force the two obstructions in Lemma~\ref{lem:reduction} to vanish,
beyond the finite-central-module class of
Theorem~\ref{thm:central-transfer}. The sought conditions should
explain elimination of these recurrences structurally, rather than
restate the vanishing of the corresponding Sylvester kernels.
\end{problem}

Theorem~\ref{thm:central-transfer} supplies this conclusion for a finite
module over a commutative central subring with the same identity,
without Noetherianity. Theorem~\ref{thm:central-separation} supplies
universal power separation when that central subring is Noetherian.
The examples and discussion in Section~\ref{sec:ring} show that either
one-sided Noetherian hypothesis alone is insufficient. Whether
two-sided Noetherianity alone suffices in the stated
Harte-quasinilpotent setting remains unresolved here.

\begin{problem}[Sharp bounds on augmented data]\label{prob:augmented-sharp}
Fix a norm and a quantitatively specified bounded class of coupled
outer-data tuples $(a,b,c,a',c')$ in a complex unital Banach algebra,
with the operators $U,V$ invertible. Determine when a finite uniform
constant $K$ exists, and its optimal value, in
\[
 \|r-p\|\le K\|pb-br\|
\]
on that class. At an exactly intertwined base, what are the sharp
first-order constants when additional constraints restrict the
feasible directions? Bounds on $a',c',U^{-1},V^{-1}$ and $b^{-1}$
are sufficient to control the coefficient in
Corollary~\ref{cor:estimate}; a sharper constrained estimate need
not attain that product bound.
\end{problem}

Theorem~\ref{thm:tangent-sharp} settles the pointwise tangent constant
for the entire local augmented class, and
Example~\ref{ex:tangent-calibration} gives its value in the Euclidean
operator norm. This does not determine the optimal uniform nonlinear
constant on a prescribed bounded class, nor the sharp constant when
its constraints exclude some of those directions.

A classical unbounded model illustrates the domain issue. On $\ell^2$,
let $A=\diag(n+1)$ and $C=\diag(-n)$, $n\ge1$, with common dense domain
\[
 D=\left\{x\in\ell^2:\sum_{n\ge1}(n+1)^2|x_n|^2<\infty\right\}.
\]
Both operators are closed and bijective with bounded inverses, and
$A+C=I|_D$, not the everywhere-defined identity. With $b=I$, all
four Hua expressions have common domain $D(A^2)$ and equal
\[
 T=A-A^2=AC=CA=C-C^2=\diag\bigl(-n(n+1)\bigr).
\]
The diagonal formula gives a closed, densely defined output with
$T^{-1}=\diag\bigl(-1/(n(n+1))\bigr)=A^{-1}+C^{-1}$.
This supplies a nonsingular model, not a singular generalized-inverse
theory.

\begin{problem}[Domain-compatible unbounded coupling]\label{prob:unbounded-coupling}
Formulate a Hua coupling theory for singular closed densely defined
operators with specified generalized-inverse data. Which domain and
compatibility hypotheses make the core restrictions and the relevant
sums and products well defined, make the Hua output closed, and yield
an inverse criterion separating its core and zero-spectral part?
Equality of the sum with $b^{-1}$ must specify its domain: if the
ordinary sum has domain the whole Banach space and both inputs are
closed, the closed graph theorem makes both inputs bounded.
The qualitative graph-gap theorem in~\cite[Theorem 4.1]{CGKR02}
motivates the continuity aspect of this question, but supplies no
Hua coupling by itself.
\end{problem}

For the last question, first recall the common-core case, which is
already settled.

Let $\A$ be a complex unital Banach algebra, and let coupled
outer data have common idempotent $e$. Put $f=1-e$ and
$T_0=et=ete$. The element $T_0$ is invertible in $e\A e$,
whose unit is $e$, with inverse $s$. For sufficiently small
nonzero $\lambda$, corner inversion gives
\[
 (T_0-\lambda1)^{-1}
 =(T_0-\lambda e)^{-1}_{e\A e}-\lambda^{-1}f.
\]
Consequently
\[
 \lim_{\lambda\to0}(et-\lambda1)^{-1}e
 =\lim_{\lambda\to0}(T_0-\lambda e)^{-1}_{e\A e}=s.
\]
Here $T_0$ is group invertible in $\A$ with group inverse $s$;
no generalized Drazin invertibility of the full output $t$ is
required. This is the common-core instance of the established
projected-resolvent and prescribed-range/kernel constructions
in~\cite[Theorems~4.1--4.2, pp.~629--631]{DjordjevicStanimirovic01}.

Distinct output supports alone also admit an elementary construction.
In a complex unital Banach algebra, let $P,Q$ be idempotents and suppose
\[
 V\in Q\A P,\qquad W\in P\A Q,\qquad WV=P,\qquad VW=Q.
\]
For $T\in\A$ with $TP=QT$, assume $H=WT$ is invertible in
$P\A P$. Then $G=H^{-1}_{P\A P}W$ satisfies
\[
 GT=P,\qquad TG=Q,\qquad GTG=G,\qquad
 G=\lim_{\lambda\to0}(WT-\lambda1)^{-1}W.
\]
Indeed, $TP=VH$ gives the product identities, and for small nonzero
$\lambda$ the complementary pole is annihilated by $W$, so
\[
 (WT-\lambda1)^{-1}W=(H-\lambda P)^{-1}_{P\A P}W.
\]
For example, $P=E_{11}$, $Q=E_{22}$, $V=T=E_{21}$ and $W=E_{12}$
in $M_2(\C)$ give $G=E_{12}$ with distinct supports. This is an outer
inverse, not a Drazin inverse of the nilpotent $T$.

For the original commuting outer input data, however,
Theorem~\ref{thm:common} shows that $ts=p$ and $st=r$, with
$s=a'+c'$, already force $p=r$, without a smallness hypothesis.
Thus weakening nilpotence or quasinilpotence alone cannot produce
distinct input cores while retaining both prescribed product identities.

\begin{problem}[Hua-compatible projected resolvents with distinct supports]
\label{prob:distinct-idempotents}
For a specified enlargement of the Hua input class, characterize when
the coupling supplies distinct output idempotents and compatible maps
$V,W$ as above, and when the resulting outer inverse $G$ equals the
designated Hua candidate. Specify which commuting-outer input
requirements, candidate formula, or prescribed product identities
are changed, together with the analytic hypotheses for the
projected-resolvent limit. The unresolved issue is this compatibility,
not the existence of an abstract distinct-support construction.
\end{problem}

\appendix
\section{Supported coupling and elementary boundaries}\label{app:supported}
The assumption that $b$ is a unit can be replaced by a sufficient confinement to its invertible Drazin corner. This appendix records that construction, not a necessity theorem for every possible generalized Hua identity.
\begin{lemma}\label{lem:transport}
Let $b\in\R^{\Dz}$, $q=bb^{\Dz}$, and suppose $a\in\R^{\Dz}$ satisfies $qa=aq=a$. Then $a^{\Dz}\in S:=q\R q$, and $a$ is Drazin invertible in the ring $S$ (with unit $q$) with the same Drazin inverse. Moreover $bq=qb\in S^\times$ with inverse $b^{\Dz}$.
\end{lemma}
\begin{proof}
From $qa=aq$ we get $q\in\comm(a)$, hence $a^{\Dz}q=qa^{\Dz}$ since $a^{\Dz}\in\comm^2(a)$~\cite{Drazin58}. Then $qa^{\Dz}=q(a^{\Dz})^2a=(a^{\Dz})^2qa=(a^{\Dz})^2a=a^{\Dz}$, and symmetrically $a^{\Dz}q=a^{\Dz}$. Consequently $a^{\Dz}\in S$. The commuting outer equations for $(a,a^{\Dz})$ hold in $S$, and the nilpotent input remainder $a-a^2a^{\Dz}$ is unchanged. These are the finite-Drazin conditions in $S$. Finally $(bq)b^{\Dz}=bb^{\Dz}q=q=b^{\Dz}(qb)$.
\end{proof}
\begin{theorem}\label{thm:supported}
Let $b\in\R^{\Dz}$, $q=bb^{\Dz}$, and let $a\in\R^{\Dz}$ satisfy $qa=aq=a$. Set $c:=b^{\Dz}-a$ and assume $c\in\R^{\Dz}$ and $aa^{\Dz}=cc^{\Dz}=:e$. Then, with $t:=a-aba$ and $s:=a^{\Dz}+c^{\Dz}$, all conclusions of Theorems~\ref{thm:common} and~\ref{thm:output} hold within the ring $S=q\R q$ upon replacing $b^{-1}$ by $b^{\Dz}$ and $1$ by $q$; explicitly, $e\le q$ commutes with $a,c,b,t$;
\[
 s=c^{\Dz}b^{\Dz}a^{\Dz}=a^{\Dz}b^{\Dz}c^{\Dz},\quad ts=st=e,\quad sts=s,\quad s=e\,(q-e+te)^{-1}_S,
\]
and $s=t^\#$ (resp.\ $t^{\Dz}$) if and only if $\nu:=(q-e)t(q-e)$ vanishes (resp.\ is nilpotent).
\end{theorem}
\begin{proof}
By Lemma~\ref{lem:transport}, $a$, and likewise $c$ (note $qc=cq=c$ since $qb^{\Dz}=b^{\Dz}q=b^{\Dz}$ and $qa=aq=a$), are Drazin invertible in $S$ with unchanged inverses, and $\bar b:=bq\in S^\times$ with $\bar b^{-1}=b^{\Dz}$, so $c=\bar b^{-1}-a$ inside $S$. Moreover $a\bar ba=a(qbq)a=(aq)b(qa)=aba$, so $t$ is unchanged. Apply Theorems~\ref{thm:common} and~\ref{thm:output} in the unital ring $S$. To compare finite-Drazin invertibility in $S$ and in $\R$, use the commuting outer equations and nilpotence of $t-t^2w$. For $t,w\in S$ these equations and this nilpotence condition are unchanged under inclusion into $\R$.

Conversely, if $w=t^{\Dz}$ in $\R$, then $w\in\comm^2(t)$ commutes with $q$. Since $qt=tq=t$ and $w=w^2t$, one obtains $qw=wq=w$, so $w\in S$ and the same equations restrict to $S$. The group case is the same argument with $t-t^2w=0$.

The indices in the imported conclusions are taken in $S=q\R q$.
More generally, for every $x\in S$ that is finite-Drazin invertible in $S$,
\[
 \ind_{\R}(x)=
 \begin{cases}
  \ind_S(x),&q=1,\\
  \max\{1,\ind_S(x)\},&q\ne1.
 \end{cases}
\]
Indeed, the inverse and all positive-exponent Drazin equations are
unchanged under inclusion, while the exponent-zero equations use $q$
and $1$, respectively. If $q\ne1$, no element supported in $q\R q$
is a unit of $\R$. In particular, when $e=q\ne1$,
$\ind_S(t)=0$ and $\ind_{\R}(t)=1$.
\end{proof}
\begin{corollary}\label{cor:supported-commuting}
Let $a,b\in\R^{\Dz}$ with $ab=ba$, and suppose $b^{\Dz}-a\in\R^{\Dz}$. Then the idempotents $aa^{\Dz}$, $bb^{\Dz}$ and $(b^{\Dz}-a)(b^{\Dz}-a)^{\Dz}$ commute pairwise; their product
\[
 e:=aa^{\Dz}\,bb^{\Dz}\,(b^{\Dz}-a)(b^{\Dz}-a)^{\Dz}
\]
is an idempotent commuting with $a,b,a^{\Dz},b^{\Dz},(b^{\Dz}-a)^{\Dz}$; the element $e\bigl(a^{\Dz}+(b^{\Dz}-a)^{\Dz}\bigr)e$ is group invertible with group idempotent $e$; and
\[
 e\,(a-aba)\,e\;=\;\Bigl[\,e\bigl(a^{\Dz}+(b^{\Dz}-a)^{\Dz}\bigr)e\,\Bigr]^\#.
\]
When $a$, $b$ and $b^{-1}-a$ are invertible, $e=1$ and this is Hua's identity.
\end{corollary}
\begin{proof}
Since $a^{\Dz}\in\comm^2(a)$ and $b\in\comm(a)$, $a^{\Dz}$ commutes with $b$; then $b^{\Dz}\in\comm^2(b)$ commutes with $a$ and $a^{\Dz}$; every element of $\{a,b,a^{\Dz},b^{\Dz}\}$ commutes with $b^{\Dz}-a$, hence with $(b^{\Dz}-a)^{\Dz}\in\comm^2(b^{\Dz}-a)$. These elements commute pairwise, so their three core idempotents commute and their product $e$ is an idempotent commuting with all of them. In the corner $e\R e$ (unit $e$): $ae$ is invertible with inverse $a^{\Dz}e$ (as $aa^{\Dz}e=e$), $be$ with $b^{\Dz}e$, and $(be)^{-1}-ae=(b^{\Dz}-a)e$ with inverse $(b^{\Dz}-a)^{\Dz}e$. Lemma~\ref{lem:hua} applied in $e\R e$ gives that $x:=e\bigl(a^{\Dz}+(b^{\Dz}-a)^{\Dz}\bigr)e$ is invertible in $e\R e$ with inverse $y:=e(a-aba)e$ (commutativity collapses $ae\,be\,ae=abae=aba\,e$). Finally, if $x,y\in e\R e$ are mutually inverse in $e\R e$, then in $\R$: $xy=yx=e$, $xyx=x$, $yxy=y$, i.e.\ $y=x^\#$ with $xx^\#=e$.
\end{proof}
\begin{example}[The naive Drazin extension fails]\label{ex:naive}
In $M_2(\C)$ let $a=E_{12}$ and $b=1$, so $c=1-a\in\R^\times$ and $t=a-a^2=a$. Here $a^{\Dz}=0$, and replacing the inverses in the inverse-of-sum identity of Lemma~\ref{lem:hua} gives
\[
 \bigl(a^{\Dz}+c^{-1}\bigr)^{\Dz}
 =\bigl((1-a)^{-1}\bigr)^{-1}=1-a\ne a=t.
\]
The input core idempotents differ: $aa^{\Dz}=0\ne1=cc^{-1}$. Theorem~\ref{thm:rigidity} shows that coincidence would be forced by intertwining for these finite-Drazin inputs. There is a separate output condition even when the cores coincide. If $aa^{\Dz}=cc^{\Dz}=e$, then Theorem~\ref{thm:output} gives $(a^{\Dz}+c^{\Dz})^{\Dz}=s^\#=te$. Thus the inverse-of-sum equality with $t$ holds exactly when the residue $\nu$ vanishes, as in Example~\ref{ex:zero} and unlike Example~\ref{ex:nilpotent}. Under coincidence the supported Bott--Duffin identities remain $ts=st=e$ and $sts=s$.
\end{example}
\begin{example}[Failure of coincidence for coupled outer data]\label{ex:noncoincident}
In $M_2(\C)$ take $a=E_{11}$ and $b^{-1}=\bigl(\begin{smallmatrix}1&1\\0&1\end{smallmatrix}\bigr)$. Then $c=\bigl(\begin{smallmatrix}0&1\\0&1\end{smallmatrix}\bigr)$ is idempotent. Choose $a'=a$ and $c'=c$, giving $p=E_{11}$ and $r=c\ne p$. Here $t=a-aba=0$ while $r\ne0$, so the prescribed product identity $st=r$ fails. Consistently, (A3) fails: $pb\ne br$.

These examples show that coincidence is a constraint on coupled outer data. Theorem~\ref{thm:common} proves that coincidence implies the supported Bott--Duffin identities. In the finite-Drazin and Banach settings of Theorem~\ref{thm:rigidity}, intertwining forces that coincidence. The inverse-of-sum equality additionally requires $\nu=0$, as explained in Example~\ref{ex:naive}.
\end{example}
\begin{remark}[Quasi-inverse form]\label{rem:quasi}
Hua's identity admits a form requiring no invertibility of the two summands. For $x,y\in\R$ with $1-xy\in\R^\times$ (equivalently $1-yx\in\R^\times$, by Jacobson's lemma), define the \emph{quasi-inverse} $x^y:=(1-xy)^{-1}x=x(1-yx)^{-1}$ (the two expressions agree because $x(1-yx)=(1-xy)x$). Then
\[
 xyx=x-(1-xy)\,x^y\,(1-yx),
\]
since $(1-xy)x^y(1-yx)=x(1-yx)=x-xyx$; and whenever $1-x^yz\in\R^\times$ one has $1-x(y+z)=(1-xy)(1-x^yz)\in\R^\times$ with the shift formula $(x^y)^z=x^{y+z}$. Under the standing condition $1-ab\in\R^\times$, specializing $x=a,y=b$ with $a,b\in\R^\times$ also makes $c=b^{-1}-a=(1-ab)b^{-1}$ a unit and recovers Lemma~\ref{lem:hua} via $a^{-1}+(b^{-1}-a)^{-1}=a^{-1}(1-ab)^{-1}$. These formulas are the associative form of the quasi-inverse identities for Jordan pairs, where $x^y=B_{x,y}^{-1}(x-Q_xy)$ for the Bergman operator $B_{x,y}$; see Loos~\cite[\S3]{Loos75}. They concern quasi-inverses under the stated unit condition; the spectral generalized-inverse questions in this paper additionally involve the core and residue conditions established above.
\end{remark}

\section*{Acknowledgments}

The author gratefully acknowledges the generous hospitality and access to library resources and other facilities provided by the Harish-Chandra Research Institute (HRI), Prayagraj (Allahabad), India, during his visits to the institute. The principal ideas underlying this work were conceived through the author's independent study and engagement with the institute's library resources during this period. The author is particularly indebted to Professor Pinaki Majumdar, whose kindness and generosity made these arrangements possible, and to whom he remains deeply grateful.

\section*{Acknowledgment of AI assistance}
The author occasionally used publicly released models of Claude (Anthropic) and ChatGPT (OpenAI) to assist him in the literature review.

\end{document}